\documentclass[a4paper, reqno, 11pt]{amsart}

\usepackage{amsfonts, amsmath, amssymb, amsthm, amscd, setspace, subfigure,mathtools,stmaryrd}
\usepackage{verbatim,mdframed,graphicx,enumitem,cancel,mathrsfs}
\usepackage{pictexwd,dcpic,slashed,fancyhdr,float,array}
\usepackage[heightrounded,textwidth=400pt,top=1in,left=1in,right=1in,bottom=1in]{geometry}
\usepackage{layout}
\usepackage[utf8]{inputenc}
\usepackage[dvipsnames]{xcolor}
\usepackage[all,cmtip]{xy}
\usepackage{xr-hyper}
\usepackage[unicode,pagebackref=true]{hyperref}
\usepackage{cite}

\numberwithin{equation}{section}

\newcommand{\geh}{g }
\newcommand\be{\begin{equation}}
\newcommand\ee{\end{equation}}
\newcommand\C{ \mathbb{C}  }
\newcommand\im{i}
\renewcommand{\a}{ \alpha }
\newcommand{\gc}{ g_{ \mathbb{C}  ^2 } }
\newcommand{\oc}{ \omega _{ \mathbb{C}  ^2 } }
\newcommand{\NN}{\ell }
\newcommand{\dd}{\delta }

\newtheorem{definition}{Definition}
\newtheorem{proposition}[definition]{Proposition}

\author[Franchetti]{Guido Franchetti}
\address[GF]{Department of Mathematical Sciences, University of Bath,
  Claverton Down, Bath BA2 7AY, England, United Kingdom}
\email{\href{mailto: gf424@bath.ac.uk}{gf424[at]bath.ac.uk}, ORCID: \href{https://orcid.org/0000-0002-1511-6204}{0000-0002-1511-6204}}
\author[Krasnov]{Kirill Krasnov}
\address[KK]{School of Mathematical Sciences, University of Nottingham, Nottingham, NG7 2RD, UK}
\email{\href{mailto: kirill.krasnov@nottingham.ac.uk}{kirill.krasnov[at]nottingham.ac.uk}, ORCID: \href{https://orcid.org/0000-0003-2800-3767}{0000-0003-2800-3767}}

\title[]{Eguchi-Hanson harmonic spinors revisited}

\begin{document}
\maketitle
\begin{abstract} We revisit the problem of determining the zero modes of the Dirac operator on the Eguchi-Hanson space. It is well known that there are no normalisable zero modes, but such zero modes do appear when the Dirac operator is twisted by a $U(1)$ connection with $L^2$ normalisable curvature. The novelty of our treatment is that we use the formalism of spin-$c$ spinors (or spinors as differential forms), which makes the required calculations simpler. In particular, to compute the Dirac operator we never need to compute the spin connection. As a result, we are able to reproduce the known normalisable zero modes of the twisted Eguchi-Hanson Dirac operator by relatively simple computations. We also collect various different descriptions of the Eguchi-Hanson space, including its construction as a hyperk\"ahler quotient of $\mathbb{C}^4$ with the flat metric. The latter illustrates the geometric origin of the connection with $L^2$ curvature used to twist the Dirac operator. To illustrate the power of the formalism developed, we generalise the results to the case of Dirac zero modes on the Ricci-flat K\"ahler manifolds obtained by applying Calabi's construction to the canonical bundle of $\mathbb{C}  P^n $.  
\end{abstract}

\maketitle
\tableofcontents

\section{Introduction}

The Dirac operator and its spectrum, in particular the spectrum of its zero modes, are of great importance and interest in both differential geometry and theoretical physics. In relation to the former, the vanishing theorem due to Lichnerowicz states that there are no harmonic spinors on compact manifolds with positive scalar curvature. When the dimension of the space of harmonic spinors is not zero, it turns out to depend on the metric, \cite{hitchin:1974}, and is not a topological invariant of the manifold. There are also beautiful results pertaining to eigenvalue estimates for the Dirac operator and the relation between this and Killing spinors, see e.g. \cite{friedrich:1997} for an authoritative exposition of all these results. In physics, harmonic spinors are important in the context of the Kaluza-Klein programme, see e.g. \cite{witten:1981}, \cite{witten:1986}, where harmonic spinors on the internal space correspond to massless particles in the physical space. For all these reasons, harmonic spinors have been studied extensively, and there is a great body of literature available on the subject. 

In the context of Riemannian geometry, one is usually interested in spinors and harmonic spinors defined on compact Riemannian manifolds. However, one can also consider non-compact spaces, in particular gravitational instantons, as is the case in this paper. In this case the relevant notion is that of normalisable, more precisely square integrable, also known as $L^2$, harmonic spinors. The non-compact case has also been studied, and appropriate index formulae have been developed, see in particular \cite{Pope:1978zx}. More recently, there has been a renewed interest in such a setup, see \cite{Jante:2013kha} and \cite{franchetti:2018}, in the context of geometric models of matter \cite{Atiyah:2011fn}.

In this paper we revisit the problem of determining the zero modes of the (twisted) Dirac operator on the Eguchi-Hanson (EH) space   \cite{eguchi:1978}. This problem has been explicitly solved, as a subcase of a more general computation, in \cite{franchetti:2018}. The novelty of our treatment in this paper is that, given that EH is a K\"ahler manifold, we can describe spinors as spin-$c$ spinors. Indeed, any almost complex manifold $M$ carries a canonical spin-$c$ structure, in which the bundle of spinors is identified with the space of differential forms $\Lambda^{0,\bullet}(M)$. The natural Dirac operator on $M$ is $\sqrt{2}(\bar{\partial}+\bar{\partial}^*)$. On a K\"ahler manifold, this operator coincides with the Dirac operator defined by the lift of the Levi-Civita connection to the spinor bundle. This fact allows us to reduce the necessary computations to elementary operations involving the action of $\bar{\partial},\bar{\partial}^*$ on the space $\Lambda^{0,\bullet}(M)$.

The operator $\sqrt{2}(\bar{\partial}+\bar{\partial}^*)$ on EH admits no normalisable zero modes, so to get an interesting problem it needs to be further twisted by a $U(1)$ connection. Up to  integer multiples, there is a unique $U(1)$ connection with square integrable curvature, and it is this connection that we use for the twisting. To elucidate the geometric meaning of this connection, we illustrate how it arises from the construction of EH as a hyperk\"ahler reduction of $\mathbb{H}^2$, which we review and further detail in this paper. One of the main results of this paper is the simple explicit expression (\ref{zm-sols}) for the EH zero modes. This should be compared with the much more involved expressions in the literature \cite{franchetti:2018}. The computations that lead to this result are also quite elementary, and take not more than a page of text, again in a favourable contrast with the other existing treatments. 

To illustrate the power of the formalism that we have developed, we study the generalisation of the problem to the case of the Calabi metric on the total space of the line bundle $\mathcal{O}(-n-1)$ over $\mathbb{CP}^n$. As it is well known, the EH metric is the simplest metric in this family and corresponds to $n=1$. We find it interesting to observe that most of the properties of EH extend to the general case, and that the latter can be handled with essentially the same techniques as those used for EH. 

Another intriguing result obtained here is that, for all values of $n$, the $U(1)$ connection with $L^2$ curvature used to twist the Dirac operator coincides with the simplest Dirac zero mode of the untwisted operator. The corresponding relation in the EH case is (\ref{twistingconn}). Of course it is only possible to appreciate this fact by making  use  of the formalism of spinors as differential forms as we do here.

The organisation of this paper is as follows. We start, in Section \ref{sec:EH}, by giving several descriptions of EH as a K\"ahler and eventually as a hyperk\"ahler manifold. We first describe it as the space obtained from the Calabi construction of a Ricci flat metric on the total space of the canonical line bundle over a K\"ahler manifold with non-zero scalar curvature. We then perform a change of coordinates that puts the base and fibre coordinates of the Calabi construction on  equal footing. This puts the EH metric in the form (\ref{EH-z12}), which is the most convenient one for spinor computations. We then discuss various harmonic 2-forms on EH, and in particular  the unique (up to scale) harmonic 2-form with $L^2$ curvature, as well as its potential. We show that this 2-form is intimately related to the exterior derivative of the (metric dual of the)  Killing vector field generating the $U(1)$ isometric action on the fibres. Finally, we describe the EH metric as the hyperk\"ahler quotient of  the flat metric on $\mathbb{H}^2={\mathbb C}^4$. We were unable to find the details of this construction, in the amount necessary to compare to our other description (\ref{EH-z12}), in the literature, so we spell them out here. We find that the $U(1)$ connection on EH with $L^2$ curvature arises naturally as the $U(1)$ connection on the total space of a circle bundle over EH obtained as the level set of the hyperk\"ahler reduction.

In Section \ref{sec:Dirac} we  determine the zero modes of the twisted Dirac operator on the EH space. We first recall why  the untwisted operator admits no normalisable zero modes and review, to the extent necessary for our purposes, results about spin-$c$ spinors and their relation to the usual Dirac spinors. We then explain how to calculate the spin-$c$ Dirac operator in practice. The key point here is that on a K\"ahler manifold, by using the spinors as differential forms approach, it is never necessary to compute derivatives of the metric. Instead,  computations only involve taking exterior derivatives and metric contractions of the relevant differential forms. This method is much more  efficient than the usual approach, where spinors are treated as  column vectors on which the $\gamma$-matrices act, and which requires the explicit  computation of the spin connection.  Finally, we solve the problem of finding the Dirac zero modes on the EH space, and obtain the explicit simple expression (\ref{zm-sols}) for the resulting harmonic spinors.

In section \ref{sec:general} we generalise  our results to the case of spinors on the total space of the canonical bundle over $\mathbb{CP}^n$. We close the paper with some concluding remarks as well as two Appendices. In the first one we review the description of $\mathbb{C}^2$ as the total space of a line bundle over $\mathbb{CP}^1$. The second one  reviews  Calabi's construction. 

\section{Eguchi-Hanson space as a K\"ahler and hyperk\"ahler manifold}
\label{sec:EH}

The Eguchi-Hanson (EH) metric is a hyperk\"ahler metric defined on a 4-manifold diffeomorphic to the cotangent bundle of $S ^2 $.   It was introduced in  \cite{eguchi:1978}, where the metric is given in bi-axial Bianchi IX form
\begin{equation}
\label{geh} 
\geh =  \left( 1 - \frac{ \kappa  }{r ^4 } \right) ^{-1}  \mathrm{d} r ^2 +  \frac{r ^2 }{4} \left( 1- \frac{\kappa  }{r ^4 } \right)  \eta _3 ^2 +  \frac{r ^2 }{4} (\eta _1 ^2 + \eta _2 ^2 ) .
\end{equation} 
Here $(\eta _i )$ are left-invariant 1-forms on $SU (2) $, see (\ref{livforms}),  
and the parameter $\kappa  $ is a positive constant. Substituting $ u ^2 = r ^2 (1- (\tfrac{\kappa }{r ^4 } ) )$ shows that the metric is regular at $r =\kappa ^{ 1/4 } $ provided that $\psi \in [0, 2 \pi )$, $\theta \in [ 0 , \pi ] $, $\phi \in [0, 2 \pi )$. The topology of a hypersurface $\Sigma _r $ of fixed $r > \kappa ^{ 1/4 }  $ is that of a circle bundle, with the circle fibre parametrised by $\psi$ and the base $S ^2 $ parametrised by $(\theta , \phi )$. If $\psi$ had its usual range of $4 \pi $, $\Sigma _r $ would be the total space $S ^3 $ of the Hopf fibration, but due to the reduced range we have instead $\Sigma _r \simeq S ^3 / \mathbb{Z}  _2  $. The level set $r = \kappa ^{ 1/4 }   $ is a 2-sphere known as a bolt \cite{gibbons:1979b}.  Asymptotically $\geh $ approaches the ($ \mathbb{Z}  _2 $-quotient of the) metric of Euclidean 4-space. Therefore, EH gives a resolution of the singularity of the $\mathbb{R}  ^4/\mathbb{Z}  ^2 = \mathbb{C}  ^2/\mathbb{Z}  ^2$ orbifold. In general, a metric with an asymptotic volume growth equal to that of $ E ^4 $ is known as asymptotically locally Euclidean, or ALE.

The form (\ref{geh}) emphasises the structure of EH as a real manifold with cohomogeneity one  under the action of $SU (2) / \mathbb{Z}  _2 $. There are several other useful descriptions of EH, which emphasise its  complex structure. We are going to review the ones that are important for our purposes. 

First, EH is a K\"ahler-Einstein manifold. One description that makes this clear 
describes EH as the total space of a complex line bundle over $\mathbb{CP}^1$, with 
complex coordinates $(w,\zeta )$  parametrising the base and fibre.
Such a construction is a special case of the more general construction in \cite{calabi:1979} of a Ricci-flat K\"ahler metric on the total space of the canonical line bundle over a K\"ahler-Einstein manifold. A second description that we need trades the ``asymmetrical'' $(w, \zeta )$ coordinates for coordinates $(z _1 , z _2 ) $ of equal standing, which play a similar role to the global  coordinates of $\mathbb{C}  ^2 $ and make the similarities between the two manifolds apparent. For comparison, the analogous descriptions of $\mathbb{C}  ^2 $ as a real cohomogeneity one manifold and as a line bundle are reviewed in Appendix \ref{c2stuff}. Finally we review (and make explicit) the description of EH as a hyperk\"ahler quotient of $\mathbb{C}  ^4 $ with its flat metric. This is of interest to us because it exhibits the $U(1)$ connection with $L^2$ curvature that we use for the twist of the Dirac operator as arising geometrically, in the total space of a certain line bundle over EH. The 5-dimensional total space in question is the level set of the hyperk\"ahler quotient construction.

\subsection{EH as a line bundle}
The paper \cite{calabi:1979}  describes two constructions leading to K\"ahler-Einstein metrics on the total space of a complex line bundle over a K\"ahler manifold. Applied to $\mathbb{C}  P ^n $, the first construction leads to a family of Ricci-flat K\"ahler metrics on the total space of the canonical line bundle $K =\mathcal{O} (-n-1 ) \rightarrow \mathbb{C}  P ^n $, the second one to a family of hyperk\"ahler metrics on $T ^\ast \mathbb{C}  P ^n $. For $n =1 $,  $\mathcal{O} (-2 )\simeq T ^\ast \mathbb{C}  P ^1 $ and both constructions result in the same metric, which is in fact the EH one.  The first construction is reviewed in Appendix \ref{calabiconstr}.  Here we  apply it to the case of the canonical bundle $K$ over $M =\mathbb{C}  P ^1  $ with the Fubini-Study metric. In this case $K$  is simply the cotangent bundle  $ \Lambda ^{ 1,0 } (\mathbb{C}  P ^1 )$. It is well known that line bundles over $\mathbb{C}  P ^1 $ are classified by their Chern number and  $K$ has Chern number $-2$.

In terms of the inhomogeneous coordinate $w $,   the  Fubini-Study metric on $\mathbb{C}  P ^1 $ takes the form
\begin{equation}
\label{cp1metric} 
g _{ \mathbb{C}  P ^1 } = \frac{4 |\mathrm{d} w |^2 }{(1 + |w |^2 )^2 },
\end{equation} 
and is isometric to the round metric on the 2-sphere of unit radius.
The corresponding K\"ahler form is
\begin{equation}
\omega _{ \mathbb{C}  P ^1 }= \frac{2i \mathrm{d} w \wedge \mathrm{d} \bar{w} }{(1 + |w| ^2 )^2 } =\frac{i}{2} e \wedge \bar e=\mathrm{vol} ,
\end{equation} 
where $\mathrm{vol}$ is the Riemannian volume element with the natural orientation as a complex manifold, and $e$ is the unitary section of $K $ given by
\begin{equation}
e =\frac{2 \mathrm{d} w }{1 + |w| ^2 }.
\end{equation} 

The metric (\ref{cp1metric}) is Einstein with scalar curvature $s =2 $, hence the Ricci form $\rho_ { \mathbb{C}  P ^1 }$, K\"ahler form $\omega _{ \mathbb{C}  P ^1 }$ and  curvature  $\mathrm{d} \a $ of the Chern connection $\a $ on $K $, see Appendix \ref{chernconn},  are related by
\begin{equation}
 \rho _ { \mathbb{C}  P ^1 }= \omega _{ \mathbb{C}  P ^1 }= -  i \mathrm{d} \a  .
\end{equation} 
It is convenient to set
\begin{equation}
\a =2 i \, a,
\end{equation} 
where
\begin{equation}
a = \frac{1}{2i}\frac{\bar{w} \mathrm{d} w - w \mathrm{d} \bar{w}   }{1 + |w| ^2 } = \operatorname{Im}  \left(\frac{\bar{w}  \mathrm{d} w }{1 + |w| ^2 }  \right) .
\end{equation} 

We now apply the construction of \cite{calabi:1979}, as described in \cite{salamon:1990}. Thus, as shown in Appendix \ref{calabiconstr}, if  $\zeta$ is a coordinate on the fibres of $K $ and we  define
\begin{align}
\theta &=\mathrm{d} \zeta + 2i \zeta a ,\\
\label{omegatotspace} 
 \omega &=2 (u \omega _{ \mathbb{C}  P ^1 }+ i u ^\prime \theta \wedge \bar \theta),
\end{align} 
where $u$ is a function of $|\zeta |^2 $ only,  then $ K  $ with K\"ahler form (\ref{omegatotspace})  is Ricci-flat K\"ahler provided that
\begin{equation}
2u u ^\prime =1 \quad \Rightarrow \quad 
u = \sqrt{  \kappa    +  |\zeta |^2 },
\end{equation} 
for $\kappa $  an integration constant.  The  associated K\"ahler form and metric are
\begin{align}
 \omega &= 2 u\,  \omega _{ \mathbb{C}  P ^1 }+i  u ^{ -1 } \theta \wedge \bar \theta ,\\
g &
=2 u \, g _{ \mathbb{C}  P ^1 }  +  2u ^{ -1 } |\theta  |^2 .
\end{align} 
It is convenient to rescale $ \zeta \rightarrow \zeta /8$, $\kappa \rightarrow \kappa /64 $, getting
\begin{align} 
\label{EHkform} 
{\omega} &= \frac{ i}{8} \left(  \left( |\zeta |^2 + \kappa  \right) ^{1/2}  e\wedge \bar{e} +  \left( | \zeta |^2 + \kappa  \right) ^{-1/2}  \theta\wedge \bar{\theta}\right)  ,\\
\label{EHmetric}
\geh &= \left( |\zeta |^2 +\kappa  \right) ^{1/2} \frac{| \mathrm{d} w|^2}{(1+|w|^2)^2} + \left( |\zeta |^2 + \kappa  \right) ^{-1/2} \frac{1}{4} | \mathrm{d} \zeta  +  2i\zeta a|^2 . 
\end{align}

Equation (\ref{EHmetric}) clearly displays EH as a non-trivial complex line bundle over $\mathbb{C}  P ^1 $ with a  twisted product metric on the total space.  The $SU (2) $ and line bundle structure of equations (\ref{geh}), (\ref{EHmetric}) should be compared with  the corresponding expressions for $\mathbb{C}  ^2 $, given by (\ref{su2metricc2}) and (\ref{gc2bundle}).
We now check that  (\ref{EHmetric})  is indeed the same metric as  (\ref{geh}). Switching to polar coordinates  $ \zeta = R e^{\im \chi }$  in the fibres gives
\be\label{EH-z-R-psi}
g= \left( R^2 + \kappa  \right) ^{1/2} \frac{| \mathrm{d} w|^2}{(1+|w|^2)^2} + \frac{1}{\left( R^2 + \kappa  \right) ^{1/2}} \frac{R^2}{4}\left( \frac{\mathrm{d} R^2}{R^2} + ( \mathrm{d} \chi +  2 a)^2 \right).
\ee
Making the coordinate change $r ^2  = R  $ shows that (\ref{EH-z-R-psi}) asymptotically becomes the flat metric on $\mathbb{C}  ^2 / \mathbb{Z}  _2 $, cfr.~(\ref{R4-zr}).  This  suggests defining
\begin{equation}
\label{raiusdef} 
r ^2 =\left(  R ^2 + \kappa    \right) ^{ 1/2 },
\end{equation} 
which gives
\begin{equation}
g =  \left( 1- \frac{ \kappa}{ r ^4} \right) ^{-1}  \mathrm{d} r^2 
+r ^2   \frac{| \mathrm{d} w|^2}{(1+|w|^2)^2} 
+ \frac{ r ^2}{4} \left( 1- \frac{ \kappa}{ r ^4} \right) ( \mathrm{d} \chi +  2 a)^2 .
\end{equation} 
Further setting
 $ w =\cot \left( \frac{\theta }{2} \right)  \mathrm{e} ^{ i \phi }$,
 \begin{equation}
\label{chiandpsi} 
 \chi =  \psi - \phi,
 \end{equation} 
 gives
\begin{equation} 
 \frac{| \mathrm{d} w|^2}{(1+|w|^2)^2}  =\frac{  \eta _1 ^2 + \eta _2 ^2 }{4}, \quad 
a =  \left( \frac{  1 + \cos \theta }{2} \right) \mathrm{d} \phi, \quad \mathrm{d} \chi + 2a = \eta _3 
\end{equation} 
so that we recover (\ref{geh}). The difference  between (\ref{chiandpsi}) and  (\ref{chiandphisu2}) is due to the different range $\psi \in [0,4 \pi )$ for  the $SU (2) $ orbits in $\mathbb{C}  ^2 $ and $\psi \in [0, 2 \pi )$ for the $SU (2) / \mathbb{Z}  _2 $ orbits in EH.

If we introduce the frame
\begin{equation}
\label{e1e2} 
e _1 = \frac{1}{2} (\kappa + |\zeta |^2 )^{ 1/4 } e, \qquad 
e _2 = \frac{1}{2} (\kappa + |\zeta |^2 )^{ -1/4 } \theta ,
\end{equation} 
the metric and K\"ahler form take the flat-space form
\begin{equation}
\label{metandkfe1e2} 
\begin{split} 
\omega &=\frac{i}{2} (e _1 \wedge \bar e _1 + e _2 \wedge \bar e _2 ),\\
\geh &= |e _1 |^2 + |e _2 | ^2 .
\end{split} 
\end{equation} 

\subsection{Equal standing coordinates}
The simplest expression for the $ \mathbb{C}  ^2 $ metric is  of course $|\mathrm{d} z _1 |^2 + |\mathrm{d} z _2 |^2  $ where the two  coordinates $(z _1 , z _2 )$ have equal standing. In the case of $\mathbb{C}  ^2 $ the bundle coordinates $(w, \zeta )$ are related to  $(z _1 , z _2 )$  by (\ref{compltobundle}). Essentially, $w =z _1 / z _2 $ is an inhomogeneous coordinate on the base of the line bundle $ \mathbb{C}  ^2 \setminus \{ 0 \} \rightarrow \mathbb{C}  P ^1 $ while $ \zeta $ parametrises the $\mathbb{C}  $ fibre.  In the case of $\mathbb{C}  ^2 $ we have $|\zeta |= \sqrt{ |z _1 |^2 + |z _2 |^2 }=r $, while for EH $ |\zeta |=R \approx r ^2 $. This suggests to define
\begin{equation} 
\label{ehtobundle} 
(z_1, z _2 ) = \frac{\sqrt{ \zeta} }{\sqrt{ 1 + |w |^2 } } (w , 1 ), \qquad 
(w, \zeta ) =\left( \frac{z _1 }{z _2 }, ( |z _1 |^2 + |z _2 |^2 ) \frac{z _2 }{\bar z _2 }\right) .
\end{equation} 
One calculates
\begin{equation}
2 \mathrm{d} z _1 =  \frac{\sqrt{\zeta }}{\sqrt{1+|w|^2}} e + \frac{w}{\sqrt{\zeta }\sqrt{1+|w|^2}} \theta, \qquad 
2 \mathrm{d} z _2 =\frac{1}{\sqrt{ \zeta }\sqrt{1+|w|^2}} \theta - \frac{\bar{w}\sqrt{\zeta }}{\sqrt{1+|w|^2}} e,
\end{equation} 
with inverse
\begin{equation}
e= \frac{2}{|\zeta |}\frac{\bar{z}_2}{ z_2}( z_2 dz_1 - z_1 dz_2), \qquad
\theta = 2\frac{z_2}{\bar{z}_2}( \bar{z}_1 dz_1 + \bar{z}_2 dz_2),
\end{equation} 
leading to
\begin{equation}
\label{EH-z12}
\geh = \frac{1}{s} \Big( F  |  z_1 \mathrm{d} z_2 - z_2 \mathrm{d} z_1 |^2+ 
F^{-1} |\bar{z}_1 \mathrm{d} z_1 + \bar{z}_2 \mathrm{d} z_2|^2 \Big),
\end{equation} 
where we have set
\begin{equation}
\label{capitalfdef} 
s= |\zeta |=|z_1|^2+|z_2|^2, \quad 
F (s)  =\sqrt{ 1 + \frac{\kappa }{s ^2 } }.
\end{equation} 

The metric (\ref{EH-z12}) should be compared with the $\mathbb{C}  ^2 $ metric in the form  (\ref{gc2split}), that is
\begin{equation}
\label{c2z1z2} 
\gc
= \frac{1}{s} \left[
 |z _1 \mathrm{d} z _2-  z _2 \mathrm{d} z _1 |^2 +   |\bar{z} _1  \mathrm{d} z _1 + \bar{z} _2 \mathrm{d} z _2 |^2 \right] ,
\end{equation} 
 which (\ref{EH-z12}) reduces to for  $\kappa =0 $. 
For later usage we note that in components (\ref{EH-z12}) becomes
\be\label{EH-expanded}
\begin{split} 
\geh_{z_1\bar{z}_1} &
=\frac{1}{s} (F |z_2|^2 + F^{-1} |z_1|^2), \qquad 
\geh_{z_2\bar{z}_2} = \frac{1}{s} (F |z_1|^2 + F^{-1} |z_2|^2), \\
 \geh_{z_1\bar{z}_2} & =  \bar \geh_{z_2 \bar{z} _1 }  =  \frac{1}{s} (F^{-1} - F) z_2 \bar{z}_1,
 \end{split} 
\ee
and, since $\geh $ has unit determinant, 
\be
\label{ehinversemetric} 
\begin{split} 
\geh^{-1}_{z_1\bar{z}_1} &= \frac{1}{s} (F |z_1|^2 + F^{-1} |z_2|^2), \qquad
 \geh^{-1}_{z_2\bar{z}_2} =  \frac{1}{s}(F |z_2|^2 + F^{-1} |z_1|^2),\\ 
  \geh^{-1}_{z_1\bar{z}_2} & =  \bar \geh ^{-1} _{z_2 \bar{z} _1 } 
  =  \frac{1}{s}   (F - F^{-1} )\bar{z}_1  z_2 .
  \end{split} 
\ee

 Writing
\begin{equation}
Z =\begin{pmatrix}
 z _1  \\
z _2 
\end{pmatrix} , \quad
J =\begin{pmatrix}
0 &1 \\
-1 & 0
\end{pmatrix} ,
\end{equation} 
we obtain
\begin{equation}
\label{gehvecform} 
\geh = \frac{1}{s} \Big( F  |  Z ^T J  \mathrm{d} Z  |^2
+  F^{-1} |Z ^T \mathrm{d} \bar Z |^2 \Big).
\end{equation} 
Clearly $s$ and $Z ^T \mathrm{d} \bar Z $ are invariant under  $Z \mapsto G Z $ for any $G \in U (2) $, while $Z ^T  J \mathrm{d} Z $ is invariant for $G \in Sp (2, \mathbb{C}  ) =SL(2, \mathbb{C}  )$. However $|Z ^T J \mathrm{d} Z |$ is  invariant for any $G \in GL(2, \mathbb{C}  ): | \det G|=1$, hence (\ref{gehvecform}) shows that  the isometry group of EH is $U(2)$.

\subsection{Harmonic 2-forms}
With respect to the $(z _1 , z _2 ) $ coordinates, the frame (\ref{e1e2}) becomes
\be\label{e12-z12}
\begin{split} 
e_1 &= \sqrt{\frac{F}{s}} \frac{\bar{z}_2}{ z_2}( z_2 dz_1 - z_1 dz_2) = \frac{1}{2} (F s )^{ 1/2 } e, \\
e_2&= \frac{ 1 }{\sqrt{F s} } \frac{z_2}{\bar{z}_2}( \bar{z}_1 dz_1 + \bar{z}_2 dz_2)=\frac{1}{2} (F s )^{ -1/2 } \theta .
\end{split} 
\ee
Since
\begin{equation}
\begin{split} 
* e _1  &= \frac{1}{2} e _1 \wedge e _2 \wedge \bar e _2, \quad * e _2  = \frac{1}{2} e _2 \wedge e _1 \wedge \bar e _1,\\
* ( e _1 \wedge \bar e _1 )&= e _2 \wedge \bar e _2 , \quad * ( e _2 \wedge \bar e _2 )= e _1 \wedge \bar e _1 , 
\end{split} 
\end{equation} 
the combination $e _1 \wedge \bar e _1 + e _2 \wedge \bar e _2 $  is self-dual. In fact, it is essentially the K\"ahler  form $\omega$ on EH, see (\ref{metandkfe1e2}).  Of course $\omega$ is also closed, hence harmonic, but not $L ^2 $.
Using (\ref{e12-z12}) $\omega$ can be rewritten in the form
\begin{equation}
\omega =2 i \left[ 
F (\mathrm{d} z _1 \wedge \mathrm{d} \bar{z} _1 + \mathrm{d} z _2 \wedge \mathrm{d} \bar{z} _2 ) 
+ \frac{1}{s}  (F ^{-1} - F ) (\bar{z} _1 \mathrm{d} z _1 + \bar{z} _2 \mathrm{d} z _2 )\wedge (z _1 \mathrm{d} \bar{z} _1 + z _2 \mathrm{d} \bar{z} _2 )
\right] .
\end{equation} 
Since $F$ satisfies
\begin{equation}
\label{Fidentities} 
F ^\prime s =F ^{-1} - F = - \frac{\kappa }{s \sqrt{ \kappa + s ^2 }} =- \frac{\kappa }{F s ^3 },
\end{equation} 
we can write the K\"ahler form as the exterior derivative of a local potential,
\begin{equation}
\label{kfpot} 
\omega = 2 i \mathrm{d} \left(  
F ( z _1 \mathrm{d} \bar{z} _1 + z _2 \mathrm{d} \bar{z} _2 ) \right).
\end{equation}

Consider now the anti self-dual combination $e _1 \wedge \bar e _1 - e _2 \wedge \bar e _2 $. It is  not closed so we look for a closed multiple of it,
\begin{equation}
\tilde \omega =\frac{i}{2} f(e _1 \wedge \bar e _1 - e _2 \wedge \bar e _2 ),
\end{equation} 
where $f $ is a function of $ s $ only. Since $\omega$ is closed, $\mathrm{d} (s F e \wedge \bar e ) = - \mathrm{d} ((s F )^{-1} \theta \wedge \bar \theta )$, so using (\ref{dzeta2}) we have
\begin{equation}
\begin{split} 
\mathrm{d} \tilde \omega &
=(2f (s F )^\prime + (sF) f ^\prime )\wedge  (\zeta \bar\theta  + \bar \zeta \theta ) \wedge e \wedge \bar e ,
\end{split} 
\end{equation} 
showing that $\tilde \omega $ is closed provided that
\begin{equation}
f =\frac{1}{(s F )^2 }.
\end{equation} 
Therefore the form
\begin{equation}
\tilde \omega =\frac{ 1}{(s F )^2 } \frac{i}{2} (e _1 \wedge \bar e _1 - e _2 \wedge \bar e _2 )
\end{equation} 
is anti self-dual and closed, hence harmonic. It is also $L ^2 $ since
\begin{equation}
| \tilde \omega | ^2 \mathrm{vol} =- \tilde \omega \wedge \tilde \omega 
= - \frac{1}{4} (s F )^{ -4 } e _1 \wedge \bar e _1 \wedge e _2 \wedge \bar e _2 \quad \Rightarrow \quad 
| \tilde \omega | ^2 = (s F )^{ -4 }
\end{equation} 
which has a finite integral over EH. Being closed, $\tilde \omega $ can also be written as the exterior derivative of a local potential, and a computation shows that
\begin{equation}
\label{asdfpot} 
 \tilde \omega 
 =  2i \mathrm{d} \left[ \frac{1}{ F s ^2 } ( z _1 \mathrm{d} \bar z _1 + z _2 \mathrm{d} \bar{z} _2   )\right] .
\end{equation} 

It is well known that the  harmonic cohomology of EH is non-trivial only in dimension two, and that the space of harmonic $L ^2 $ 2-forms is 1-dimensional \cite{hitchin:2000}. As we have just seen  it is generated by $\tilde \omega $.
% One can also check that the top holomorphic form $\Lambda $ becomes
%\be\label{Omega}
%4 dz_2\wedge dz_1 = \theta \wedge e = \Omega.
%\ee

It is interesting to compare the expressions of $\omega$, $\tilde \omega $ with that of $\mathrm{d} \theta _3 $, for $\theta _3 $ the metric dual with respect to the EH metric of the Killing vector field $X _3 $,
\begin{equation} 
\label{adaptedcoframeeh} 
\begin{split}
%\theta _1 &= X _1 ^\flat 
%=+ \frac{ i}{4} \left(  z _1 \mathrm{d} {z} _2 - z _2 \mathrm{d} {z} _1  -\bar{z}  _1  \mathrm{d} \bar z _2  + \bar{z}  _2 \mathrm{d} \bar z _1 \right) 
%= - \frac{1}{2} \operatorname{Im} (z _1 \mathrm{d} z _2 - z _2 \mathrm{d} z _1 ),\\
%\theta _2 &= X _2 ^\flat 
%=+ \frac{ 1}{4} \left(  z _1 \mathrm{d} {z} _2 - z _2 \mathrm{d} {z} _1  + \bar{z}  _1  \mathrm{d} \bar z _2 -\bar{z}  _2 \mathrm{d} \bar z _1 \right) 
%= + \frac{1}{2} \operatorname{Re} (z _1 \mathrm{d} z _2 - z _2 \mathrm{d} z _1 ),\\
\theta_3  &=X_3 ^\flat  
= \frac{ i}{2F} \left(  z _1 \mathrm{d} \bar{z} _1  + z _2 \mathrm{d} \bar{z} _2  -\bar{z}  _1  \mathrm{d} z _1 -\bar{z}  _2 \mathrm{d} z _2 \right)
=  \frac{1}{F} \operatorname{Im} (\bar{z} _1 \mathrm{d} z _1 + \bar{z} _2 \mathrm{d} _2 ) .
%,\\
%\theta _4  &=X _4 ^\flat  
%=- \frac{1}{4}  \left(  z _1 \mathrm{d} \bar{z} _1  + z _2 \mathrm{d} \bar{z} _2   + \bar{z}  _1  \mathrm{d} z _1  + \bar{z}  _2 \mathrm{d} z _2 \right) 
%= -\frac{1}{2}  \operatorname{Re}  (\bar{z} _1 \mathrm{d} z _1 + \bar{z} _2 \mathrm{d} _2 ) .
\end{split} 
\end{equation} 
Since $X _3 $ is a Killing vector field and EH is Ricci-flat, $\mathrm{d} \theta _3  $ is  harmonic.  One calculates
\begin{equation}
\label{dtheta3} 
\begin{split} 
2\mathrm{d} \theta _3  &
= \frac{i}{2} \left[ e _1 \wedge \bar e _1 + e _2 \wedge \bar e _2 - \kappa  \left( \frac{ e _1 \wedge \bar e _1 - e _2 \wedge \bar e _2}{ (s F )^2  }  \right)  \right] 
= \omega  - \kappa \tilde \omega ,
\end{split} 
\end{equation} 
hence $\kappa \tilde \omega $ represents the same cohomology class as $\omega$, which generates of $H ^2 _{ \mathrm{dR}  }(\mathrm{EH}  )$. We can also see that $\mathrm{d} \theta _3 $ is harmonic but not $L ^2 $ and, interestingly, that $\omega$, $\kappa \tilde \omega $ are the self-dual and anti self-dual parts of $\mathrm{d} \theta _3 $,
\begin{equation}
\begin{split}
\omega &= * \mathrm{d} \theta _3 + \mathrm{d} \theta _3  , \qquad 
\kappa \tilde \omega =* \mathrm{d} \theta _3  - \mathrm{d} \theta _3 .
\end{split}
\end{equation} 
Using (\ref{kfpot})  and (\ref{asdfpot}) we see that $\mathrm{d} \theta _3 $ can also be written in the form
\begin{equation}
\mathrm{d} \theta _3 
= i \mathrm{d} \left[  \frac{1}{F} (z _1 \mathrm{d} \bar z _1 + z _2 \mathrm{d} \bar{z} _2) \right] .
\end{equation} 

\subsection{$U (1) $ connection with $L ^2 $ harmonic curvature} In Section \ref{harmspinorseh} we will consider the Dirac operator on EH twisted by a $U (1) $ connection $\mathcal{A} $ with $L ^2 $ harmonic curvature $\mathrm{d} \mathcal{A}  $. As we just discussed, $\mathrm{d} \mathcal{A}  $ is necessarily some constant multiple of $\tilde \omega $ and by (\ref{asdfpot})  we can take  $\mathcal{A} $ to be some multiple of
\begin{equation} 
 \frac{1 }{s ^2 F }  ( z _1 \mathrm{d} \bar{z} _1 + z _2 \mathrm{d} \bar{z} _2  - \bar{z} _1 \mathrm{d} z _1 - \bar{z} _2 \mathrm{d} z _2 ).
\end{equation} 
In order for $\mathrm{d} \mathcal{A} $ to be the curvature of a connection we need to impose the quantisation condition
\begin{equation}
\frac{i}{2 \pi } \int _{ \mathbb{C}  P ^1 }    \mathrm{d} \mathcal{A}   =\NN  \in \mathbb{Z} ,
\end{equation} 
obtaining
\begin{equation}
\label{twistingconn} 
\mathcal{A} =\NN  (A- \bar A ), \quad 
A = \frac{ \sqrt{ \kappa } }{2s ^2 F } (z _1 \mathrm{d} \bar{z} _1 + z _2 \mathrm{d} \bar{z} _2  )= \frac{1}{2}\,  \partial  \sinh ^{-1}  (\sqrt{ \kappa }/s  ) .
\end{equation}

\subsection{Hyperk\"ahler quotient}
\label{subsec:hyperkquot}
It is well known, see e.g. \cite{hitchin:1992}, that EH can be obtained as the hyperk\"ahler quotient of $\mathbb{H}  ^2 $, but we were unable to find anywhere in the literature the details of this construction in the amount sufficient for comparison with the two descriptions given above. We spell it out here.\footnote{One of us (KK) benefited from a discussion with Daniel Platt in relation to the material described in this subsection.}

The first part of this construction is standard and appears in many references. Identify $\mathbb{H} ^2 $ with $ \mathbb{C}  ^4  $ via the isomorphism $\mathbb{H}  = \mathbb{C} ^2  +  \mathbb{C}  ^2 \mathbf{j}$, $q _i = Z _i +   W_i \,  \mathbf{j}$, $i =1, 2 $, $q _i \in \mathbb{H}  $, $Z _i , W _i \in \mathbb{C}  $, and equip $\mathbb{H} ^2  $ with the flat metric and the K\"ahler forms
\begin{equation}
\begin{split} 
\omega _R &=\frac{i}{2} (\mathrm{d} Z _1 \wedge  \mathrm{d} \bar{Z} _1 + \mathrm{d} Z _2 \wedge  \mathrm{d} \bar{Z} _2 + \mathrm{d} W _1 \wedge  \mathrm{d} \bar{W} _1 + \mathrm{d} W _2 \wedge  \mathrm{d} \bar{W} _2 ) =\frac{i}{2}\partial \overline{\partial } (|Z |^2 + |W |^2 ),\\
\omega _ C &
=\omega _2 + i \omega _3 = \mathrm{d} W _1 \wedge \mathrm{d} Z _1  + \mathrm{d} W _2 \wedge \mathrm{d} Z _2 ,
\end{split} 
\end{equation} 
where $Z =(Z _1 , Z _2 )^T $, $W =(W _1 , W _2 ) ^T $, $|Z |^2 =|Z _1 |^2 + |Z _2 |^2 $.

The right $U (1) $ action
\begin{equation} 
\label{rightu1c4} 
(Z  , W )\mapsto (Z, W )\mathrm{e} ^{ \mathbf{i} t } =(\mathrm{e} ^{ \mathbf{i} t }  Z , \mathrm{e} ^{ -\mathbf{i} t }W ),
\end{equation} 
which corresponds to translation along the $U (1) $ fibres of the Hopf fibration $S ^1 \hookrightarrow S ^{ 7 }  \rightarrow \mathbb{C}  P ^3 $, is Hamiltonian and isometric. 
%It is generated by the Killing vector field
%\begin{equation}
%X =\frac{i}{2} \left(  z _1 \frac{\partial }{\partial z _1 } + z _2 \frac{\partial }{\partial z _2 } - w _1 \frac{\partial }{\partial w _1 } - w _2 \frac{\partial }{\partial w _2 } - cc\right) ,
%\end{equation} 
%where $ cc $ stands for complex conjugate.
% The associated moment maps, determined by $i _X \omega  =\mathrm{d} \mu  $, are (up to constant factors I am not being careful about)
It is convenient to take the associated moment maps to be
\begin{equation}
\mu _R = | Z |^2 - |W |^2 -2 \sqrt{ \kappa }, \quad \mu _{C}  =Z ^T W =Z _1 W _1 + Z _2 W _2 .
\end{equation} 
The  level set $\mu ^{-1} (0) =\mu _R ^{-1} (0) \cap \mu _C ^{-1} (0) $ is a smooth real 5-manifold.  Further quotienting by the $U (1) $ action (\ref{rightu1c4}) we get the hyperk\"ahler quotient
\begin{equation}
\mathbb{H}  ^2 \sslash U (1) =\mu ^{-1} (0) / U (1) .
\end{equation} 

To compare to our previous description (\ref{EH-z12}) of  EH, we need to parametrise this level set by  coordinates related to $z_1, z_2$, as well as some coordinate for the $S^1$ fibre. This is something we were unable to find in the literature. To proceed, let $h \in U (2) $ act on $\mathbb{C}  ^2 $ by ordinary matrix multiplication. Then  the isometric left $U (2)$ action on $\mathbb{C}  ^4 $ given by
\begin{equation}
\label{u2leftaction}
 h \cdot (Z,W ) = (h Z, \bar h W )
\end{equation} 
 commutes with (\ref{rightu1c4}) and preserves $\mu$  level sets, so it descends to  an isometric action on the quotient. In quaternionic notation, the moment map is
\be\label{moment-map}
\mu (q_1, q_2) = \frac{1}{2} \sum_{a=1,2} q_a {\bf i} \bar{q}_a - \sqrt{ \kappa  }\, \mathbf{i} 
\ee
and
(\ref{u2leftaction}) corresponds to 
\begin{equation}
h \cdot (q _1 , q _2 )  = (q _1 , q _2 )h ^T 
\end{equation} 
where on the rhs we have ordinary matrix multiplication.

The isometric $U (2 )$ action on the quotient just introduced should match the one described in the paragraph following \eqref{gehvecform}. This is   achieved by setting
\begin{alignat}{2}
\label{parametr-wz}
Z_1 &= z_1 e^{\im\psi} f^{1/2},& \qquad 
Z_2 &= z_2 e^{\im\psi} f^{1/2}, \\ \nonumber 
 W_1&= - z_2 e^{-\im\psi} f^{-1/2},& \qquad
 W_2 &= z_1 e^{-\im\psi} f^{-1/2},
\end{alignat}
where $f$ is some function of $z_1, z_2$ to be determined. This ansatz automatically satisfies $ \mu _C  (Z,W) = 0 $, while the condition $|Z |^2 - |W |^2 = 2 \sqrt{  \kappa } $  becomes
\be
(|z_1|^2+|z_2|^2) (f-f^{-1}) = 2\sqrt{\kappa},
\ee
which is solved by
\be
\label{ffun} 
f = \sqrt{1+ \frac{\kappa}{s^2}} + \frac{\sqrt{\kappa}}{s}.
\ee

We note the useful identities
\be
f- f^{-1} = \frac{2\sqrt{\kappa}}{s}, \qquad f+ f^{-1} = 2\sqrt{1+ \frac{\kappa}{s^2}} = 2F,
\ee
as well as
\be
\frac{1}{f^{1/2}} \frac{\mathrm{d} f^{1/2}}{\mathrm{d}s} = \frac{1}{2f} \frac{\mathrm{d}f}{\mathrm{d}s} = - \frac{\sqrt{\kappa}}{2s^2 F}.
\ee

We now pull-back the flat metric on $\C^4$ to $\mu^{-1}({\bf i}\sqrt{\kappa})$ as parametrised by (\ref{parametr-wz}). We have
\be
\begin{split} 
\mathrm{d}Z_1 &= z_1 e^{\im\psi} f^{1/2}\left( \frac{\mathrm{d}z _1}{z_1} + \im \mathrm{d}\psi    - \frac{\sqrt{\kappa}}{2s^2 F} \mathrm{d}s   \right),  \quad 
\mathrm{d}Z _2  = z_2 e^{\im\psi} f^{1/2}\left( \frac{\mathrm{d}z _2}{z_2} + \im \mathrm{d}\psi    - \frac{\sqrt{\kappa}}{2s^2 F} \mathrm{d}s   \right),
\\ 
\mathrm{d}W _1 &= -z_2 e^{-\im\psi} f^{-1/2}\left( \frac{\mathrm{d}z _2}{z_2} - \im \mathrm{d}\psi    + \frac{\sqrt{\kappa}}{2s^2 F} \mathrm{d}s   \right), \quad 
\mathrm{d}W _2 = z_1 e^{-\im\psi} f^{-1/2}\left( \frac{\mathrm{d}z _1}{z_1} - \im \mathrm{d}\psi    + \frac{\sqrt{\kappa}}{2s^2 F} \mathrm{d}s   \right),
\end{split} 
\ee
giving
\be
\begin{split} 
|\mathrm{d}Z _1|^2+|\mathrm{d}Z _2 |^2 &
= f (|\mathrm{d}z _1|^2+|\mathrm{d}z _2|^2) + s f \mathrm{d}\psi  ^2 + \im f \mathrm{d}\psi   (z_1 \mathrm{d}\bar{z}_1- \bar{z}_1 \mathrm{d}z _1 + z_2 \mathrm{d}\bar{z}_2- \bar{z}_2 \mathrm{d}z _2) \\  &
+ f \frac{\kappa (\mathrm{d}s   )^2}{4s^3 F^2} - f \frac{\sqrt{\kappa}(\mathrm{d}s   )^2}{2s^2 F}, \\ 
|\mathrm{d}W _1 |^2+|\mathrm{d} W _2 |^2 &= f^{-1} (|\mathrm{d}z _1|^2+|\mathrm{d}z _2|^2) + s f^{-1} \mathrm{d}\psi  ^2 - \im f^{-1} \mathrm{d}\psi   (z_1 \mathrm{d}\bar{z}_1- \bar{z}_1 \mathrm{d}z _1 + z_2 \mathrm{d}\bar{z}_2- \bar{z}_2 \mathrm{d}z _2) \\ &
+ f^{-1} \frac{\kappa (\mathrm{d}s   )^2}{4s^3 F^2} + f^{-1} \frac{\sqrt{\kappa}(\mathrm{d}s   )^2}{2s^2 F}.
\end{split} 
\ee
This means that  the pull-back of (half) the flat metric on $\C^4$ is 
\be
\begin{split} 
&\frac{1}{2}( |\mathrm{d}Z _1|^2+|\mathrm{d}Z _2|^2+ |\mathrm{d}W _1|^2+|\mathrm{d}W _2|^2) 
= F (|\mathrm{d}z _1|^2+|\mathrm{d}z _2|^2)+ sF \mathrm{d}\psi  ^2 
\\ &
+ \im \mathrm{d}\psi   \frac{\sqrt{\kappa}}{s} (z_1 \mathrm{d}\bar{z}_1- \bar{z}_1 \mathrm{d}z _1 + z_2 \mathrm{d}\bar{z}_2- \bar{z}_2 \mathrm{d}z _2) 
- \frac{\kappa (\mathrm{d}s   )^2}{4s^3 F}.
\end{split} 
\ee
We now complete the square getting 
\be
\frac{1}{2}( |\mathrm{d}Z _1|^2+|\mathrm{d}Z _2|^2+ |\mathrm{d}W _1|^2+|\mathrm{d}W _2|^2) 
= \geh + sF( \mathrm{d}\psi   + \im (A-\bar{A}))^2 ,
\ee
where 
\begin{align} \label{EH-HK}
\geh &= F ( |\mathrm{d}z _1|^2 + |\mathrm{d} _2|^2) + \frac{1}{s}(F^{-1}-F)|\bar{z}_1 \mathrm{d}z _1 + \bar{z}_2 \mathrm{d}z _2|^2,\\
\label{apotential} 
A&= \frac{\sqrt{\kappa}}{2 s^2 F} (z_1 \mathrm{d}\bar{z}_1+z_2 \mathrm{d}\bar{z}_2).
\end{align} 
Using the identity
\begin{equation}
\label{identityc2}
s \Big ( |\mathrm{d} z _1 |^2 + |\mathrm{d} z _2 |^2 \Big ) = |z _1 \mathrm{d} z _2 - z _2 \mathrm{d} z _1 |^2 + |\bar{z} _1 \mathrm{d} z _1 + \bar{z} _2 \mathrm{d} z _2 |^2 
\end{equation} 
  we see that (\ref{EH-HK})  is the EH metric in the form (\ref{EH-z12}). 
We recognise  $A$ for the $(0,1)$ part of the ${\rm U}(1)$ connection  (\ref{twistingconn}), whose curvature is anti-self-dual and $L ^2 $. Therefore,  the natural connection $ A-\bar{A}$ on EH  arises as the connection in the 5-dimensional fibered space $\mu^{-1}({\bf i}\sqrt{\kappa})$  obtained in the process of a hyperk\"ahler reduction of $\mathbb{H}  ^2$.

\section{Dirac zero modes  on the Eguchi-Hanson space}
\label{sec:Dirac}

\subsection{No Dirac zero modes on EH}

\label{harmspinorseh} 
Let  $D$ be the Dirac operator on EH and $\sigma$ be a spinor. Since the scalar curvature $s $ of EH vanishes, Lichernowicz's identity
\begin{equation} 
D ^2 \sigma = \nabla ^\ast \nabla  \sigma + \frac{s}{4} \sigma
\end{equation} 
implies that EH admits no $L ^2  $ Dirac zero modes. In order to obtain a non-trivial problem it is necessary to twist the spinor bundle by a complex line bundle equipped with a connection $\mathcal{A} $. Equivalently, we need to replace the spin structure on EH by a spin-$c $ structure. Lichernowicz's identity then becomes 
\begin{equation}
\label{lichner-twisted}
D _\mathcal{A}  ^2 \sigma = \nabla ^\ast \nabla  \sigma + \frac{s}{4} \sigma +  \mathrm{d} \mathcal{A}   \cdot \sigma,
\end{equation} 
where $\cdot $ denotes Clifford multiplication, and for a suitable choice of $\mathcal{A} $ it is possible to obtain non-trivial solutions of the twisted Dirac equation $D _\mathcal{A}  \sigma =0 $.
Twisting by an arbitrary connection does not make for an interesting problem, but taking $\mathcal{A} $ so that $\mathrm{d} \mathcal{A} $ is  $L ^2 $ harmonic is a natural choice. As we already discussed, in the case of EH any such connection $\mathcal{A} $ takes the form (\ref{twistingconn}). Therefore, we want to find $L ^2 $ solutions of the equation $D _ \mathcal{A}  \sigma =0 $ for $\mathcal{A} $ given by (\ref{twistingconn}). Working on the spinor bundle, viewed as rank four Spin(4) module,  this problem has been considered in \cite{franchetti:2018}. 

\subsection{The canonical spin-$c$ structure}

However, since EH is a K\"ahler manifold, there is a more convenient approach  in terms of complex differential forms.  In fact, any almost complex manifold carries a canonical spin-$ c $ structure. Let us recall the necessary background. We follow
\cite{friedrich:1997}, section 3.4.

We have the following proposition:
\begin{proposition} The spinor bundle $S$ of an Hermitian manifold of complex dimension $k$ (with respect to an arbitrary spin-$c$ structure) is isomorphic to
$$ S= (\Lambda^{0,0}\oplus \ldots \oplus \Lambda^{0,k})\otimes S_0= (\Lambda^{0,0}\oplus \ldots \oplus \Lambda^{k,0})\otimes S_k,$$
where
$$ S_0=\{ \sigma\in S: \omega\sigma = \im k \sigma\}, \quad S_k=\{ \sigma\in S: \omega\sigma = -\im k \sigma\},$$
and $\omega$ is the K\"ahler form, acting on a spinor by Clifford multiplication. In particular $S_0 = \Lambda^{k,0}\otimes S_k$ and 
$S_k = \Lambda^{0,k}\otimes S_0$.
\end{proposition}

The {\it anti-canonical spin-$c$ structure} on an Hermitian manifold $M,J$ is the one for which the bundle $S_0$ is trivial and $S_k$ coincides with the canonical bundle $K=\Lambda^k(T^*)$ of $M$. We have the following proposition, see \cite{friedrich:1997} and also \cite{nicolaescu:2007}.
\begin{proposition}
    Let $(M,J)$ be a K\"ahler manifold equipped with the anti-canonical spin-$c$ structure. Then 
    $$ S\simeq \Lambda^{0,0}\oplus \ldots \oplus \Lambda^{0,k},$$
    and the Dirac operator defined by the Levi-Civita connection coincides with 
    \begin{equation} 
\sqrt{ 2 } (\overline{\partial } + \overline{\partial }^\ast ).
\end{equation} 
\end{proposition}

For reasons that will become clear in the next subsections,  we will consider instead the Dirac operator
\begin{equation} 
\label{ourdiracope} 
D =\overline{\partial } - \overline{\partial }^\ast ,
\end{equation} 
which has the same kernel as $\overline{\partial } + \overline{\partial }^\ast$, with the Clifford action given by
\begin{equation}
\label{clactionforms} 
\upsilon  \cdot \sigma = \upsilon ^{0,1 } \wedge \sigma  + \iota _{ \upsilon  ^\sharp } \sigma  .
\end{equation}

We discussed above how the spin structure on EH admits no $L ^2  $ harmonic spinors. The same holds for the (anti-)canonical spin-$c$ structure. In fact the spinor bundle is isomorphic to $ W = \Lambda ^0 \oplus \Lambda ^{ 0,1 } \oplus \Lambda ^{ 0,2 }$ and harmonic spinors correspond to forms $\sigma  \in W $ such that 
\begin{equation}
D \sigma  =(\overline{\partial }- \overline{\partial }^\ast )\sigma  =0 .
\end{equation}
On a complete K\"ahler manifold
\begin{equation}
D \sigma  =0 \Leftrightarrow D ^2 \sigma  =0 \Leftrightarrow  (\mathrm{d} \mathrm{d} ^\ast + \mathrm{d} ^\ast \mathrm{d} ) \sigma  =0 \Leftrightarrow \mathrm{d} \sigma  =0 =\mathrm{d} ^\ast \sigma  .
\end{equation} 
Thus harmonic spinors in $\Lambda ^0 $ are constant functions and  $\alpha \in \Lambda ^2 $ is harmonic if and only if $* \alpha $ is a constant function. Neither is $L ^2 $. As for 1-forms, by the Bochner identity
\begin{equation}
\triangle \sigma   =\nabla ^\ast \nabla \sigma  + \operatorname{Ric} (\sigma  ).
\end{equation} 
Since EH is Ricci-flat, any harmonic form has to be parallel, hence cannot be $L ^2 $. The result is not surprising  since the Chern connection on the canonical bundle of a Ricci-flat manifold has zero curvature, so the Dirac operator associated to the canonical spin-$c$ structure on EH is equivalent to the one associated to the spin structure.

\subsection{Computing the action of the Dirac operator on spin-$c$ spinors}

In four dimensions, the two chiralities of spin-$c$ spinors are identified with
\be
W^+ = \Lambda^0 + \Lambda^{0,2}, \qquad W^- = \Lambda^{0,1},
\ee
and the spin-$c$ Dirac operator is
\be
D = \bar{\partial} - \bar{\partial}^*.
\ee
Here $\bar{\partial}$ is the $(0,1)$ projection of the exterior derivative $\mathrm{d}$, and $\bar{\partial}^*$  its adjoint. The use of the minus sign here rather than plus is a convention that is more suitable for our purposes, because it leads to a very simple  way of computing $D$, which we will explain below.

The action of $\bar{\partial}$ on $\Lambda^0, \Lambda^{0,1}$ is very simple
\begin{eqnarray}
\bar{\partial} : \Lambda^0 \ni \alpha \to \partial_{\bar{z}_1} \alpha \, d\bar{z}_1 + \partial_{\bar{z}_2} \alpha \,d\bar{z}_2 \in \Lambda^{0,1}, \\ \nonumber
\bar{\partial} : \Lambda^{0,1} \ni \beta \, d\bar{z}_1 + \gamma \, d\bar{z}_2 \to (\partial_{\bar{z}_1} \gamma -  \partial_{\bar{z}_2} \beta) d\bar{z}_1 \wedge d\bar{z}_2 \in \Lambda^{0,2}.
\end{eqnarray}
Here $\alpha,\beta,\gamma$ are functions of $z_1,\bar{z}_1,z_2,\bar{z}_2$. The action of $\bar{\partial}^*$ involves the metric and is more complicated. However, on a K\"ahler manifold there is a very simple  set of rules to be followed. We derive these rules after $\bar{\partial}^*$ is computed. 

\subsection{Computation of $\bar{\partial}^*$}

We first compute the action of $\bar{\partial}^*$ on $\Lambda^{0,1}$. For the sake of generality, we do the computation for a general Hermitian metric on $\C^2$. The Hermitian pairing of $\bar{\partial}\alpha \in \Lambda^{0,1}$ with an arbitrary element $\beta \, d\bar{z}_1 + \gamma \, d\bar{z}_2 \in \Lambda^{0,1}$ is given by
\begin{eqnarray}
\langle \bar{\partial} \alpha, \beta \, d\bar{z}_1 + \gamma \, d\bar{z}_2 \rangle = g^{-1}( \partial_{\bar{z}_1} \alpha \, d\bar{z}_1 + \partial_{\bar{z}_2} \alpha \,d\bar{z}_2, \bar{\beta} \, dz_1 + \bar{\gamma} \, dz_2) = \\ \nonumber
 \partial_{\bar{z}_1} \alpha \, \bar{\beta} g^{-1}_{z_1\bar{z}_1} + \partial_{\bar{z}_1} \alpha \, \bar{\gamma}g^{-1}_{z_2\bar{z}_1} +
 \partial_{\bar{z}_2} \alpha \, \bar{\beta} g^{-1}_{z_1\bar{z}_2} 
 + \partial_{\bar{z}_2} \alpha \, \bar{\gamma} g^{-1}_{z_2\bar{z}_2}.
 \end{eqnarray}
 We now multiply this by the volume factor, which is $v_g |dz_1|^2 |dz_2|^2$, and integrate by parts. This gives
 \be\nonumber
v_g \langle \alpha, \bar{\partial}^* (\beta \, d\bar{z}_1 + \gamma \, d\bar{z}_2)\rangle =
-\alpha \overline{\partial_{z_1}\left( g^{-1}_{z_1\bar{z}_1} v_g \beta \right)} -\alpha \overline{\partial_{z_1}\left(  g^{-1}_{z_1\bar{z}_2}  v_g \gamma\right) }
-\alpha\overline{\partial_{z_2}\left(  g^{-1}_{z_2\bar{z}_1} v_g \beta \right)}- \alpha\overline{\partial_{z_2}\left( g^{-1}_{z_2\bar{z}_2} \gamma \right)},
\ee
which is an identity that only holds modulo surface terms (which we assume to vanish). 
We can simplify this expression using the fact that the components of the inverse metric times the volume form give the components of the original metric. Doing so we get
\be
v_g \bar{\partial}^* (\beta \, d\bar{z}_1 + \gamma \, d\bar{z}_2) = -\partial_{z_1}\left( g_{z_2\bar{z}_2}  \beta \right) +\partial_{z_1}\left(  g_{z_2\bar{z}_1}  \gamma\right) 
+\partial_{z_2}\left(  g_{z_1\bar{z}_2} \beta \right)- \partial_{z_2}\left( g_{z_1\bar{z}_1} \gamma \right)\in \Lambda^0.
\ee

The computation of $\bar{\partial}^*$ on $\Lambda^{0,2}$ is similar. We have
\begin{eqnarray}
\langle (\partial_{\bar{z}_1} \gamma -  \partial_{\bar{z}_2} \beta) d\bar{z}_1 \wedge d\bar{z}_2, \delta d\bar{z}_1 \wedge d\bar{z}_2\rangle = \\ \nonumber
g^{-1}( ( \partial_{\bar{z}_1} \gamma -  \partial_{\bar{z}_2} \beta) (d\bar{z}_1 \otimes d\bar{z}_2- d\bar{z}_2 \otimes d\bar{z}_1), \bar{\delta} (dz_1 \otimes dz_2)) = \\ \nonumber
 ( \partial_{\bar{z}_1} \gamma -  \partial_{\bar{z}_2} \beta)v_g^{-1} \bar{\delta}.
 \end{eqnarray}
 We now multiply by $v_g$ and integrate by parts to get
 \be
 v_g \langle \beta d\bar{z}_1 + \gamma d\bar{z}_2, \bar{\partial}^* ( \delta d\bar{z}_1 \wedge d\bar{z}_2) \rangle
=-  \gamma \overline{ \partial_{z_1} \delta} +  \beta \overline{ \partial_{z_2} \delta}.
\ee
Since we want to rewrite $\bar{\partial}^* ( \delta d\bar{z}_1 \wedge d\bar{z}_2) $ in the form
\be
v_g \bar{\partial}^* ( \delta d\bar{z}_1 \wedge d\bar{z}_2) = A d\bar{z}_1+Bd\bar{z}_2,
\ee
we need to match
\be
-  \gamma \overline{ \partial_{z_1} \delta} + \beta \overline{ \partial_{z_2} \delta}=
\beta \bar{A} g^{-1}_{z_1\bar{z}_1} + \beta \bar{B} g^{-1}_{z_2\bar{z}_1} + \gamma\bar{A} g^{-1}_{z_1\bar{z}_2} + \gamma\bar{B} g^{-1}_{z_2\bar{z}_2},
\ee  
which gives
\be
A = g_{z_1\bar{z}_1} \partial_{z_2}  \delta - g_{z_2\bar{z}_1} \partial_{z_1}  \delta, \qquad
B = g_{z_1\bar{z}_2} \partial_{z_2}  \delta - g_{z_2\bar{z}_2} \partial_{z_1} \delta.
\ee

Collecting the above results, we obtain the  two chiral halves of the Dirac operator. Abusing the notation, we  still denote by $D:W^+\to W^-$ one of the two chiral parts, and  
by $D^\dagger: W^-\to W^+$ the other one. Their action is
\begin{eqnarray}
\label{D-plus} D(\alpha + \delta d\bar{z}_1\wedge d\bar{z}_2) =  
\left( \partial_{\bar{z}_1}\alpha -v_g^{-1} g_{z_1\bar{z}_1} \partial_{z_2} \delta + v_g^{-1}g_{z_2\bar{z}_1} \partial_{z_1} \delta\right) d\bar{z}_1 \\ \nonumber
+\left(  \partial_{\bar{z}_2}\alpha - v_g^{-1}g_{z_1\bar{z}_2} \partial_{z_2} \delta + v_g^{-1}g_{z_2\bar{z}_2} \partial_{z_1} \delta \right) d\bar{z}_2,
\end{eqnarray}
\begin{eqnarray}\label{D-minus}
D^\dagger( \beta d\bar{z}_1+\gamma d\bar{z}_2) = 
(\partial_{\bar{z}_1} \gamma -  \partial_{\bar{z}_2} \beta) d\bar{z}_1 \wedge d\bar{z}_2 \\ \nonumber
+v_g^{-1} \partial_{z_1}\left( g_{z_2\bar{z}_2} \beta \right) - v_g^{-1}\partial_{z_1}\left(  g_{z_2\bar{z}_1}  \gamma\right) 
-v_g^{-1}\partial_{z_2}\left(  g_{z_1\bar{z}_2} \beta \right)+ v_g^{-1}\partial_{z_2}\left( g_{z_1\bar{z}_1} \gamma \right) .
\end{eqnarray}

\subsection{Spin-{\it c} Dirac operator on a K\"ahler manifold}
So far we have computed the Dirac operator only assuming that the metric is Hermitian. If the metric  is also K\"ahler, 
thanks to the  identities
\be
\partial_{z_1} g_{z_2 \bar{z}_1} = \partial_{z_2} g_{z_1 \bar{z}_1}, \qquad
\partial_{z_1} g_{z_2 \bar{z}_2} = \partial_{z_2} g_{z_1 \bar{z}_2},
\ee
the terms involving derivatives of the metric in (\ref{D-minus}) cancel among each other.
As a result, we obtain the following simplified expression for the chiral Dirac operator on $W^-$,
\begin{eqnarray}\label{D-minus*}
D^\dagger( \beta d\bar{z}_1+\gamma d\bar{z}_2) = 
(\partial_{\bar{z}_1} \gamma -  \partial_{\bar{z}_2} \beta) d\bar{z}_1 \wedge d\bar{z}_2 \\ \nonumber
+v_g^{-1} g_{z_2\bar{z}_2} \partial_{z_1}\beta  - v_g^{-1}g_{z_2\bar{z}_1} \partial_{z_1}\gamma
-v_g^{-1}g_{z_1\bar{z}_2}  \partial_{z_2} \beta + v_g^{-1}g_{z_1\bar{z}_1} \partial_{z_2} \gamma .
\end{eqnarray}

\subsection{Practical way of computing the Dirac operator on a K\"ahler manifold}

The computations above illustrate  that there is a very simple method for computing the Dirac operator  on a K\"ahler manifold. The method consists in taking any differential form in question, and then first computing its {\bf full} exterior derivative, including the terms not belonging to $\Lambda^{0,k}$. One then uses the metric pairing to map the latter terms from $\Lambda^{1,k}$ to $\Lambda^{0,k-1}$. Let us see how this works in practice. 

We start with the spinor in $W^-=\Lambda^{0,1}$. We have
\begin{eqnarray}
\mathrm{d} ( \beta \mathrm{d} \bar{z}_1+\gamma \mathrm{d} \bar{z}_2) = \partial_{z_1} \beta \mathrm{d} z_1\wedge \mathrm{d} \bar{z}_1 + \partial_{z_2} \beta \mathrm{d} z_2\wedge \mathrm{d} \bar{z}_1 + \partial_{z_1} \gamma \mathrm{d} z_1\wedge \mathrm{d} \bar{z}_2 + \partial_{z_2} \gamma \mathrm{d} z_2\wedge \mathrm{d} \bar{z}_2 \\ \nonumber
+ (\partial_{\bar{z}_2} \gamma - \partial_{\bar{z}_1} \beta) \mathrm{d} \bar{z}_1\wedge \mathrm{d} \bar{z}_2.
\end{eqnarray}
The second line here is in $\Lambda^{0,2}$, and so a spinor in $W^+$, but
the first line lies in $\Lambda^{1,1}$ and is not a spinor.  However, we can use the metric pairing to map it into $\Lambda^0$. Denoting this projection by $g^{-1}$ we have
\begin{eqnarray}
g^{-1}\left(\mathrm{d}( \beta \mathrm{d}\bar{z}_1+\gamma \mathrm{d}\bar{z}_2) \right) = g^{-1}_{z_1\bar{z}_1} \partial_{z_1} \beta
+g^{-1}_{z_2\bar{z}_1} \partial_{z_2} \beta + g^{-1}_{z_1\bar{z}_2} \partial_{z_1} \gamma+g^{-1}_{z_2\bar{z}_2} \partial_{z_2} \gamma
\\ \nonumber
+ (\partial_{\bar{z}_2} \gamma - \partial_{\bar{z}_1} \beta) \mathrm{d}\bar{z}_1\wedge \mathrm{d}\bar{z}_2.
\end{eqnarray}
Taking into account the relation between the metric and its inverse we can write this as
\begin{eqnarray}
g^{-1}\left(\mathrm{d}( \beta \mathrm{d}\bar{z}_1+\gamma \mathrm{d}\bar{z}_2) \right) = v^{-1}_g g_{z_2\bar{z}_2} \partial_{z_1} \beta
- v_g^{-1} g_{z_1\bar{z}_2} \partial_{z_2} \beta -v_g^{-1}  g_{z_2\bar{z}_1} \partial_{z_1} \gamma+ v_g^{-1} g_{z_1\bar{z}_1} \partial_{z_2} \gamma
\\ \nonumber
+ (\partial_{\bar{z}_2} \gamma - \partial_{\bar{z}_1} \beta) \mathrm{d}\bar{z}_1\wedge \mathrm{d}\bar{z}_2,
\end{eqnarray}
which is the correct expression (\ref{D-minus*}) for the chiral Dirac operator on $W^-$. 

Completely analogous computations give the other chiral Dirac operator: we first apply the full exterior derivative to $\Lambda^{0,2}$
\be
\mathrm{d}( \delta \mathrm{d}\bar{z}_1 \wedge \mathrm{d}\bar{z}_2) = \partial_{z_1} \delta \mathrm{d}z_1 \wedge \mathrm{d}\bar{z}_1 \wedge \mathrm{d}\bar{z}_2+  \partial_{z_2} \delta \mathrm{d}z_2 \wedge \mathrm{d}\bar{z}_1 \wedge \mathrm{d}\bar{z}_2.
\ee
We now do all possible metric pairings to map this into $\Lambda^{0,1}$. We have
\be
g^{-1}\left( \mathrm{d}( \delta \mathrm{d}\bar{z}_1 \wedge \mathrm{d}\bar{z}_2)\right) = g^{-1}_{z_1\bar{z}_1} \partial_{z_1} \delta\, \mathrm{d}\bar{z}_2 - g^{-1}_{z_1\bar{z}_2} \partial_{z_1}  \delta\, \mathrm{d}\bar{z}_1 
+  g^{-1}_{z_2\bar{z}_1} \partial_{z_2} \delta \, \mathrm{d}\bar{z}_2
- g^{-1}_{z_2\bar{z}_2} \partial_{z_2} \delta \, \mathrm{d}\bar{z}_1.
\ee
Again using the relation between the metric and its inverse and collecting  terms we have
\be\nonumber
g^{-1}\left( \mathrm{d}( \delta \mathrm{d}\bar{z}_1 \wedge \mathrm{d}\bar{z}_2)\right) =v_g^{-1} (g_{z_2\bar{z}_1} \partial_{z_1}  \delta -v_g^{-1} g_{z_1\bar{z}_1} \partial_{z_2} \delta )\mathrm{d}\bar{z}_1 + v_g^{-1} (g_{z_2\bar{z}_2} \partial_{z_1} \delta\, -  v_g^{-1} g_{z_1\bar{z}_2} \partial_{z_2} \delta) \mathrm{d}\bar{z}_2,
\ee
which is exactly the $\delta$-dependent part of (\ref{D-plus}). All in all, this is a very simple way of computing the spin-$c$ Dirac operator for a K\"ahler metric, which involves nothing more complicated than taking the exterior derivative, and then doing metric contractions. It is in order to have such a simple recipe for computing $D$ that we have taken the minus sign in our definition $D  =  \overline{\partial } - \overline{\partial }^\ast$  of the Dirac operator. 

To summarise, the spinor $D \sigma $ can be obtained by  first calculating    $\mathrm{d} \sigma  \in \Lambda ^{ 0, q + 1 } \oplus \Lambda ^{ 1,q }$. The $\Lambda ^{ 0, q + 1 }  $ component corresponds to  $\overline{\partial } \sigma $ while   contracting with the (inverse) metric maps the $(1, q )$ component to $ - \overline{\partial }^\ast \sigma \in \Lambda ^{ 0,q-1 }$. 

 Similarly, when considering the twisted operator
\begin{equation}
D _ \mathcal{A}  = D + \mathcal{A} , 
\end{equation} 
 the action of $\mathcal{A} $ is given by (\ref{clactionforms}) so we first calculate $ \mathcal{A}  \wedge \sigma  \in \Lambda ^{ 1,q }\oplus \Lambda ^{ 0, q+ 1 }$ and then contract the $(1,q ) $ part with the inverse metric to obtain a form of degree $(0,q-1 )$. This gives very simple computational rules, allowing to compute the twisted Dirac operator with minimal effort. In particular, we never need to compute the spin connection for EH. Nor do we ever need to compute the derivatives of the metric components. This is, of course, part of the magic of K\"ahler geometry. 

 \subsection{Zero modes of the twisted Dirac operator on EH}
 
  We now proceed with the calculation using the EH metric in the form (\ref{EH-z12}) and the connection (\ref{twistingconn}).  Write 
\begin{equation}
 W ^+ = \mathbb{C}  \oplus \Lambda ^{ 0,2 }, \quad  W ^- = \Lambda ^{ 0,1 }  
 \end{equation} 
for the even and odd part of $\Lambda ^{ 0, \bullet }$.
A generic spinor  $\sigma = \sigma _+ + \sigma _- $, $\sigma _\pm \in W ^\pm $, has the form
\begin{equation}
\sigma _+ =\alpha  + \delta\,  \mathrm{d} \bar{z} _1 \wedge \mathrm{d} \bar{z} _2 , \quad \sigma _- =\beta \, \mathrm{d} \bar{z}  _1 + \gamma \, \mathrm{d} \bar{z} _2 ,
\end{equation} 
with $\alpha  $, $\beta $, $\gamma$, $\delta$ functions of $z _i $, $\bar z _i $, $i =1,2 $.
However, we are twisting by a connection whose curvature is anti-self-dual, see (\ref{asdfpot}). Such anti-self-dual curvature can only Clifford act non-trivially in (\ref{lichner-twisted}) on spinors of one chirality. It can be checked that this are the spinors in $W^-$. This means that there are no non-trivial $L ^2 $ harmonic spinors in $ W ^+ $. Hence we take
\begin{equation}
\sigma =\sigma _- =\beta \, \mathrm{d} \bar{z}  _1 + \gamma \, \mathrm{d} \bar{z} _2.
\end{equation} 
We calculate
\begin{equation} 
\begin{split} 
\mathrm{d} \sigma  _- &= \partial_{z_1} \beta \mathrm{d}z_1\wedge \mathrm{d}\bar{z}_1 + \partial_{z_2} \beta \mathrm{d}z_2\wedge \mathrm{d}\bar{z}_1 + \partial_{z_1} \gamma \mathrm{d}z_1\wedge \mathrm{d}\bar{z}_2 + \partial_{z_2} \gamma \mathrm{d}z_2\wedge \mathrm{d}\bar{z}_2 \\ &
+ (\partial_{\bar{z}_1} \gamma - \partial_{\bar{z}_2} \beta) \mathrm{d}\bar{z}_1\wedge \mathrm{d}\bar{z}_2.
\end{split}
\end{equation}  
The second line belongs to  $W^+$ while, as discussed,  we need to contract the first line with the inverse metric, getting
\be
 D  \sigma _- 
= g^{-1}_{z_1\bar{z}_1} \partial_{z_1} \beta
+g^{-1}_{z_2\bar{z}_1} \partial_{z_2} \beta + g^{-1}_{z_1\bar{z}_2} \partial_{z_1} \gamma
+g^{-1}_{z_2\bar{z}_2} \partial_{z_2} \gamma
\\ \nonumber
+ (\partial_{\bar{z}_1} \gamma - \partial_{\bar{z}_2} \beta) \mathrm{d} \bar{z}_1\wedge \mathrm{d} \bar{z}_2.
\ee
Substituting  (\ref{ehinversemetric}) we obtain
\be
\begin{split} 
D \sigma _- &
= \frac{1}{s} (F |z_1|^2 + F^{-1} |z_2|^2) \partial_{z_1} \beta
+    \frac{1}{s}   (F - F^{-1} )  (\bar{z}_1  z_2 \partial_{z_2} \beta + \bar{z}_2  z_1  \partial_{z_1} \gamma)
+ \frac{1}{s}(F |z_2|^2 + F^{-1} |z_1|^2) \partial _{ z _2 }\gamma
\\ &
+ (\partial_{\bar{z}_1} \gamma - \partial_{\bar{z}_2} \beta)  \mathrm{d} \bar{z}_1\wedge \mathrm{d} \bar{z}_2,
\end{split} 
\ee
which can be rewritten in the form
\be\label{Dm-action-psim}
D\sigma _- = \left( \frac{F-F^{-1}}{s}\right)  ( z_1 \partial_{z_1} + z_2 \partial_{z_2}) (\bar{z}_1 \beta+\bar{z}_2 \gamma) +  F^{-1} ( \partial_{z_1} \beta + \partial_{z_2} \gamma)  \\ \nonumber
+( \partial_{\bar{z}_1} \gamma -    \partial_{\bar{z}_2} \beta) \mathrm{d} \bar{z}_1\wedge \mathrm{d} \bar{z}_2.
\ee
Consider now the action of the connection (\ref{twistingconn}). We calculate
\be
\begin{split} 
%(A_n-\bar{A}_n)\cdot \sigma _+ &
%= \frac{\alpha n\sqrt{\kappa}}{Fs^2} ( z_1 \mathrm{d} \bar{z}_1+z_2 \mathrm{d} \bar{z}_2) + \frac{2\delta n\sqrt{\kappa}}{s^2} ( \bar{z}_2 \mathrm{d} \bar{z}_1-  \bar{z}_1 \mathrm{d} \bar{z}_2),\\
\mathcal{A} \cdot \sigma _- &
 = \frac{\NN\sqrt{\kappa}}{2Fs^2}( z_1 \gamma-z_2 \beta) \mathrm{d} \bar{z}_1\wedge \mathrm{d} \bar{z}_2 - \frac{\NN\sqrt{\kappa}}{2s^2} ( \bar{z}_1 \beta + \bar{z}_2 \gamma).
 \end{split} 
 \ee
 Putting all together, the twisted Dirac equation is
 \begin{equation}
 \label{eqs-ZM} 
 \begin{split} 
 D _\mathcal{A}   \sigma _- &
 =  \left(\left( \frac{F-F^{-1}}{s}\right)  \left( z_1 \partial_{z_1} + z_2 \partial_{z_2}\right)  - \frac{\NN\sqrt{\kappa}}{2s^2}\right)  (\bar{z}_1 \beta+\bar{z}_2 \gamma)  +  F^{-1} \left(  \partial_{z_1} \beta + \partial_{z_2} \gamma \right) \\ &
 +\left( \partial _{ \bar{z} _1 } \gamma - \partial _{ \bar{z} _2 } \beta +   \frac{\NN\sqrt{\kappa}}{2Fs^2}( z_1 \gamma-z_2 \beta) \right) 
 \mathrm{d} \bar{z}_1\wedge \mathrm{d} \bar{z}_2=0.
 \end{split} 
 \end{equation}

We take the ansatz, which solves the $\Lambda^{0,2}$ part of (\ref{eqs-ZM}),
\be\label{ansatz-zero}
\beta = z_1^{N-m+1} z_2^{N+m} h(s), \qquad \gamma= z_1^{N-m} z_2^{N+m+1} h(s)
\ee
for $h(s)$  an  function of $s$ to be determined below. Here $N  \geq 0 $ is such that $2 N \in \mathbb{Z}  $,  $m =\{ - N , -N + 1 , \ldots , N \} $. The reason for this particular ansatz is that, since the left  $SU (2) $ action on $(z _1 , z _2 )$  is an isometry, spinors can be decomposed into irreducible $SU (2) $ representations. By taking $\beta$, $\gamma$ as in (\ref{ansatz-zero}) we have
\begin{equation}
\sigma _- = h  z_1^{N-m} z_2^{N+m}     (\bar{z} _1 \mathrm{d} z _1 + \bar{z} _2 \mathrm{d} z _2 ),
\end{equation} 
where $h (\bar{z} _1 \mathrm{d} z _1 + \bar{z} _2 \mathrm{d} z _2 ) $ is $SU (2) $-invariant and the space of homogeneous polynomials in $z _1 , z _2 $  of degree $2N$ gives the $SU (2) $ irrep of dimension $2 N + 1 $.

 The  $\Lambda^0$ part of (\ref{eqs-ZM}) now reduces to
\be
F^2( 2N h + (hs)') +h = \frac{\NN}{2}\sqrt{\kappa} \, \frac{hF}{s},
\ee
which is solved by
\be
h = \left( \frac{1}{Fs^{2N+2}} \right) \frac{1}{f ^{ \NN }},
\ee
where 
\be
f = \sqrt{1+ \frac{\kappa}{s^2}} + \frac{\sqrt{\kappa}}{s}
\ee
is the same function as (\ref{ffun}). Thus, we have found the following zero modes
\be\label{zm-sols}
\sigma _-  = \frac{z_1^{N-m} z_2^{N+m}}{F s^{2N+2} f ^{ \NN }}
( z_1 \mathrm{d} \bar{z}_1+z_2 \mathrm{d} \bar{z}_2).
\ee

\subsection{Normalisability}
 Let us discuss normalisability. Using (\ref{ehinversemetric}) we calculate 
 \be
 \begin{split} 
|z_1 \mathrm{d} \bar{z}_1+z_2 \mathrm{d} \bar{z}_2|^2  =
  |z_1|^2 \frac{F |z_1|^2+ F^{-1} |z_2|^2}{s} + |z_2|^2 \frac{F |z_2|^2+ F^{-1} |z_1|^2}{s} + 2\frac{F-F^{-1}}{s} |z_1|^2 |z_2|^2   = F s,
  \end{split} 
 \ee
hence
 \be\label{norm-squared}
 |\sigma _- |^2 = \frac{|z_1|^{2(N-m)} |z_2|^{2(N+m)}}{F s^{4N+3} f ^{ 2\NN }} .
 \ee
 As discussed previously,  EH asymptotically approaches the flat metric on $\mathbb{C}  ^2 / \mathbb{Z}  _2 $ with radial coordinate $ r = \sqrt{ s } = \sqrt{ |z _1 |^2 + |z _2 |^2  }$ and volume element $\sim r ^3 \mathrm{d} r $.  We need to check the $L ^2 $ condition for small and large $r$. For large $r$,  $f \sim F \sim 1 $ hence
 \begin{equation}
 |\sigma _- |^2 r ^3 \sim \frac{1}{r ^{ 4 N + 3 } }
 \end{equation} 
and the $L ^2 $ condition gives $N > -1/2 $, which is satisfied by any non-negative half-integer $N$.
 For small $r $, $f \sim  F \sim r ^{-2} $, hence
 \begin{equation}
 |\sigma _- |r ^3 \sim r ^{ 2\NN-4 N -1 },
 \end{equation} 
hence a normalisable spinor needs to satisfy  $2N< \NN $.  

In conclusion,  harmonic spinors belong to  $SU (2) $ representations of dimension $2N + 1 $. A spinor is normalisable if and only if 
 \begin{equation} 
 \label{ehl2cond} 
 1 \leq 2N + 1 \leq \NN 
 \end{equation} 
where $\NN \in \mathbb{Z}  $ is the flux of the curvature 2-form $ \mathrm{d} \mathcal{A} $.  Here we have assumed $\NN >0 $,  the case $\NN<0 $ can be treated similarly.  As we already knew, in the untwisted case $\NN=0$ there are no normalisable zero modes. For $\NN=1$ we have only the singlet zero mode, for $\NN=2$ we have both the singlet and doublet, and so on, for a total of $\NN(\NN+1)/2$ zero modes for general $\NN$.

 \section{The general case}
\label{sec:general}
 
 \subsection{Calabi's metric on $ \mathcal{O}(-n-1)$}
We now apply Calabi's construction to the canonical bundle $K =\mathcal{O} (- n - 1 )\rightarrow \mathbb{C}  P ^n $ over $\mathbb{C}  P ^n $ equipped with the Fubini-Study metric. Applying the construction of Appendix \ref{calabiconstr} to $M =\mathbb{C}  P ^n $ we get, for $\lambda \neq 0$, $\kappa >0 $  arbitrary  constants,
\begin{align}
 \omega &=\lambda u\,  \omega_{ \mathbb{C}  P ^n } +i  (n + 1 )(\lambda u )^{ -n } \theta \wedge \bar \theta ,\\
 \label{calabicpn}
g &
= \lambda u \, g _{ \mathbb{C}  P ^n }  +2 (n + 1 )  (\lambda u )^{ -n } |\theta  |^2 .
\end{align} 
Here
\begin{equation}
u =  \left(  c | \zeta |^{ 2 }  +  \kappa  \right) ^{ \frac{1}{n + 1} }, 
\end{equation} 
and, for $\zeta$  a complex coordinate on the fibres,
\begin{equation}
\theta =\mathrm{d} \zeta + \alpha \zeta 
\end{equation} 
with $\alpha$ the Chern connection on $K$. Its curvature  $\mathrm{d} \alpha $ satisfies
\begin{equation}
-i \mathrm{d} \alpha  =\frac{s _{ \mathbb{C}  P ^n }}{2n}  \omega _{ \mathbb{C}  P ^n },
\end{equation} 
where $s_{ \mathbb{C}  P ^n }$ is the scalar curvature of the FS metric. The constants $c$ and $\lambda$ are related by
\begin{equation} 
c =\frac{s_{ \mathbb{C}  P ^n } (n + 1 )^2 }{2n \lambda ^{ n + 1 }}.
\end{equation} 

Let $(w _i ) $, $i =1 , \ldots ,n $, be inhomogeneous coordinates on $\mathbb{C}  P ^n $. With respect to the local K\"ahler potential 
\begin{equation}
\mathcal{K} = \frac{C }{2} \log (1 + |w  |^2  ),
\end{equation} 
where $|w |^2 =|w _1 |^2 + \cdots + |w _n |^2 $, $C\in \mathbb{R}  ^\times $ is some constant, the FS metric and K\"ahler form take the form $(g _{ \mathbb{C}  P ^n }) _{ \mu \bar\nu } = \partial _\mu \partial _{ \bar\nu } \mathcal{K} $, 
$(\omega _{ \mathbb{C}  P ^n }) _{ \mu \bar\nu } =  i (g _{ \mathbb{C}  P ^n }) _{ \mu \bar\nu }$. The scalar curvature is related to $C$ by
\begin{equation}
\label{scurvcpn} 
 s  _{ \mathbb{C}  P ^n }= \frac{ 4  n (n + 1 ) }{C}.
 \end{equation} 
Since
 \begin{equation}
 \partial \mathcal{K} = \frac{C }{2} \frac{\bar{w} \, \mathrm{d} w }{1 + |w| ^2 }, \quad 
\overline{\partial } \mathcal{K} = \overline{\partial \mathcal{K} },
 \end{equation} 
the K\"ahler form of $\mathbb{C}  P ^n $ can be written 
\begin{equation}
\omega _{ \mathbb{C}  P ^n } 
=  \frac{i}{2} \left( \partial \overline{\partial } \mathcal{K} - \overline{\partial } \partial \mathcal{K} \right)
 =\frac{i}{2}  \mathrm{d} \left( \overline{\partial }\mathcal{K} - \partial \mathcal{K} \right) ,
\end{equation} 
so that
\begin{equation}
\begin{split}
\mathrm{d} \alpha &
= -\left(  \frac{n + 1}{C }\right) \mathrm{d}  (\overline{\partial }\mathcal{K} - \partial \mathcal{K}  ).
\end{split}
\end{equation} 
Hence up to the addition of a closed form
\begin{equation}
\label{alphaform} 
\alpha 
=2i \left( \frac{n + 1}{2} \right) \frac{\operatorname{Im} (\bar{w} \mathrm{d} w )}{1 + |w| ^2 }.
\end{equation} 
Note that $\alpha$ is purely imaginary. As we did for $n =1 $ we write
\begin{equation}
\label{varconnzz} 
\alpha =2i \, a , \quad
 a =\frac{1}{2i}   \left( \frac{n + 1}{2} \right) \left(  \frac{\bar{w} \mathrm{d} w - w \mathrm{d} \bar{w} }{1 + |w| ^2 }\right) .
\end{equation} 

\subsection{Rewriting in terms of "symmetrical" coordinates}

Introducing the homogeneous coordinates $(z_i)$, $i=1,\ldots,n+1$, related to $(z_i)$ by $w_i=\frac{z_i}{z_{n+1}}$,  the FS metric on $\mathbb{C}  P ^n $ can also be written as (the pullback along a holomorphic section of)
\begin{equation}
g _{ \mathbb{C}  P ^n } = \frac{1}{s ^2 } \sum _{ 1 \leq i<j \leq n + 1 } |z _i \mathrm{d} z _j - z _j \mathrm{d} z _i |^2 ,
\end{equation} 
where 
\begin{equation}
s = |z _1 |^2 + \cdots + |z _{n + 1 } |^2 .
\end{equation} 
We will also need the analogous of the identity (\ref{identityc2}), which for general $n$ reads 
\begin{equation}
\label{identitycn}
s  |\mathrm{d} z |^2 = | \bar{z}  \mathrm{d} z |^2  
+\!\!\!\! \sum _{ 1 \leq i<j \leq  n + 1  } | z _{ i  }\mathrm{d} z _j - z _j \mathrm{d} z _{ i  } |^2  .
\end{equation} 
We are using the notation $z=(z_1,\ldots z_{n+1})$ and denoting the dot product by juxtaposition so that e.g.
\begin{equation}
    \bar z \mathrm{d}z = 
    \bar z _1 \mathrm{d}z_1 +\cdots+\bar z _{n+1} \mathrm{d}z_{n+1} = \partial s.
\end{equation}

Having collected all the ingredients,  we would like to rewrite the Calabi metric in term of the complex coordinates $ (z _i )$. 
The  base-fibre coordinates $(w_i,\zeta)$ are related to the more ``symmetrical'' coordinates $(z_i)$ by
\begin{equation}
\label{relavzetas} 
w _i = \frac{z _i }{z _{ n + 1 }} ,\  i =1 , \ldots ,n , \qquad 
\zeta =s ^{ \frac{n + 1}{2 }}  \left( \frac{z _{ n + 1 }}{\bar z _{ n + 1 }} \right) ^{ \frac{n + 1}{2 }},
\end{equation} 
with inverse
\begin{equation}
z _i =w _i  \frac{\zeta ^{  \frac{1}{ n + 1} }}{\sqrt{ 1 + |w |^2 }}, \ i =1 , \ldots , n , \qquad 
z _{ n + 1 } =\frac{\zeta ^{  \frac{1}{ n + 1} }}{\sqrt{ 1 + |w |^2 }}.
\end{equation} 

One calculates
\begin{equation}
\label{reldiffzetas} 
\bar z \mathrm{d} z 
= |\zeta |^{ \frac{2}{n + 1} }  \frac{1}{2} \left( \frac{\bar{w} \mathrm{d} w - w \mathrm{d} \bar{w} }{1 + |w |^2 } \right)  + \left( \frac{1}{n + 1} \right) \frac{\bar\zeta \mathrm{d} \zeta }{|\zeta |^{ \frac{2n}{n + 1} }},
\end{equation} 
so that, using (\ref{varconnzz}), we have
\begin{equation}
\label{relzbardzeta} 
 (n + 1 )\bar z \mathrm{d} z 
=\frac{\bar\zeta }{|\zeta |^{ \frac{2n}{n + 1 }}} (\mathrm{d} \zeta + \alpha \zeta ) 
=\frac{\bar\zeta }{|\zeta |^{ \frac{2n}{n + 1 }}} \theta .
\end{equation} 
It follows 
\begin{equation} 
\label{thetasquared} 
|\theta |^2 
=(n + 1 )^2 s  ^{ n -1 }  |\bar z \mathrm{d} z |^2 .
\end{equation} 
%Note that
%\begin{equation}
%\frac{\bar\zeta }{|\zeta |^{ \frac{2n}{n + 1} }}
%= \frac{1}{s ^{ \frac{ n -1}{2} }} \left( \frac{ \bar Z _{ n + 1 }}{Z _{ n + 1 }} \right) ^{ \frac{n + 1}{2} }.
%\end{equation} 
Note that, since $\alpha$ is purely imaginary, taking the real and imaginary part of (\ref{relzbardzeta})  we get
\begin{equation} 
 \mathrm{d} s
= \mathrm{d} (|\zeta |^2 ) ^{ \frac{1}{1 + n} },
\end{equation} 
in agreement with (\ref{relavzetas}), and
\begin{equation} 
(n + 1 ) \operatorname{Im} (\bar z \mathrm{d} z  )
= |\zeta |^{ \frac{2}{n + 1} } \alpha   + \frac{ \operatorname{Im} (\bar\zeta \mathrm{d} \zeta) }{|\zeta |^{ \frac{2n}{n + 1} }}.
\end{equation}

We can now rewrite \ref{calabicpn} in terms of the coordinates $(z_i)$. Using (\ref{identitycn}) we have
\begin{equation} 
\begin{split} 
g &
= \lambda u \, g _{ \mathbb{C}  P ^n }  +2 (n + 1 ) (\lambda u )^{ -n } |\theta  |^2 
= \frac{ 2\lambda u }{s ^2 }\!\!\!\!\!\! \sum _{ 1 \leq i<j \leq n + 1 } |z _i \mathrm{d} z _j - z _j \mathrm{d} z _i |^2 
+ 2 (n + 1 )(\lambda u )^{ -n }  |\theta  |^2\\ &
= \frac{ 2 \lambda u}{ s ^2 } \left( s |\mathrm{d} z |^2  - |\bar z \mathrm{d} z |^2  \right) 
+2 (n + 1 ) (\lambda u )^{ -n } |\theta |^2  .
\end{split} 
\end{equation} 
We now rescale $\zeta \rightarrow \frac{ \zeta }{(n + 1)^{ 3/2 } } $, so that by (\ref{thetasquared}),
\begin{equation}
\label{rescaledtheta2} 
  |\theta |^2 \rightarrow \frac{|\theta |^2 }{(n + 1 )^3 } = \frac{s ^{ n -1 }}{n + 1} |\bar z \mathrm{d} z |^2 ,
 \end{equation} 
 hence
\begin{equation} 
g= 2 \left[ \frac{  \lambda u}{ s }  |\mathrm{d} z |^2
+ \left( \frac{ s ^{ n -1 }}{(\lambda u )^{ n} }    -\frac{  \lambda u}{ s ^2 }  \right)  |\bar z \mathrm{d} z |^2  
  \right].
\end{equation} 
We are going to drop the overall factor of $2$. Using the relation $|\zeta |=s ^{ \frac{n + 1}{2} }$  to rewrite $u$ as a function of $s$ we obtain
\begin{equation}
\frac{\lambda u }{s}  
= \lambda c ^{ \frac{1}{n + 1} } \left( 1 + \frac{\kappa  }{cs ^{n + 1 } } \right) ^{ \tfrac{1}{n + 1} }.
\end{equation} 
Rescaling $\kappa  \rightarrow c \kappa $ and choosing $c = \lambda ^{ - (n + 1 )} $ 
we define
\begin{equation}
F (s) = \frac{\lambda u}{s} = \left( 1 + \frac{\kappa }{s ^{ n + 1 }} \right) ^{ \frac{1}{ n + 1 }},
\end{equation} 
so that 
\begin{equation} 
F ^\prime =\frac{1- F ^{n + 1 } }{F ^{n} s} = \frac{s ^{ n-1 }}{(\lambda u )^n } - \frac{\lambda u }{s ^2 }.
\end{equation} 
Therefore, the Calabi metric on $\mathcal{O} (-n-1 )$   is given in terms of the coordinates $(z _i )$ by
\begin{equation}
\label{calabimetriccpn} 
g =F |\mathrm{d} z |^2  + F ^{ \prime  } \,|\bar z \mathrm{d} z |^2.
\end{equation} 
%in agreement with the calculation in \cite{lye:2022}. 
Since $F $ only depends on $s$, it is clear that the  metric is invariant under the left action of $U (n + 1) $ on $ \mathbb{C}  ^{n + 1 } $.
The  K\"ahler form corresponding to (\ref{calabimetriccpn}) is

\begin{equation}
\omega =2i \left( F \mathrm{d} z \wedge \mathrm{d} \bar z+ F ^\prime \bar z \mathrm{d} z \wedge z \mathrm{d} \bar z  \right)     =2 i \mathrm{d} \left( F  z \mathrm{d} \bar z \right) .
\end{equation} 
For small $s$, $F \sim s $ and  $ \operatorname{Re} (F z \mathrm{d} \bar z)\sim \frac{1}{2} \mathrm{d} \log s  $, so $F  z \mathrm{d} \bar z $  is not well-defined for $s =0 $  and $ \omega $ is not exact. Note that for $n =1 $  (\ref{calabimetriccpn})  reduces to
\begin{equation}
g =  F |\mathrm{d} z |^2 + \frac{1}{s} \left( \frac{1}{F } - F \right)  |\bar z \mathrm{d} z |^2 , \quad F = \sqrt{ 1  + \frac{\kappa }{s ^2 } }
\end{equation} 
and we recover (\ref{EH-HK}).

The function $F $ satisfies the identity
\begin{equation}
\label{frel} 
F ^{ \prime } = \frac{1}{s} \left( \frac{1 }{F ^{n} }  -F \right) 
\quad \Leftrightarrow \quad F ^{ n } (s F )^\prime =1.
\end{equation} 
In components
\begin{equation}
\label{metrcsec} 
g _{ \mu \bar\nu } 
= \partial _\mu (F  z_\nu  )
=F ^{  \prime }\bar z_{ \mu }  z _\nu   +F \delta _{ \mu \nu },
\end{equation} 
with determinant
\begin{equation}
\det (g _{ \mu \bar\nu } )
=F ^{ n  } (F s ) ^\prime =1,
\end{equation} 
as it should be since the metric is Ricci-flat. The inverse of the matrix $g _{ \mu \bar\nu }$  has components
\begin{equation}
\label{ginverse} 
g ^{  \rho   \bar  \nu  }
= \frac{\delta ^{ \rho \nu }}{F } + \frac{\kappa  }{s ^{ n  + 2 } F  }\bar z ^\rho  z ^\nu  .
\end{equation} 

\subsection{$U(1)$ connection with $L^2$ harmonic curvature}
The canonical spin-$c$ structure on $K$ admits no non-trivial $L ^2 $ harmonic spinor, hence to get non-trivial zero modes we need to twist by a line bundle equipped with a $U(1)$ connection $\mathcal{A}$. The curvature $\mathrm{d} \mathcal{A} $ is a purely imaginary form of degree $(1,1) $ which we  also want to be $L ^2 $ harmonic. 

The space of $L ^2 $ harmonic forms on $K$ is 1-dimensional.  In fact, in the setting of \cite{hausel:2004}, the Calabi metric on $K$ is a scattering metric, with $K$  viewed as a fibration with trivial fibre and $X = \overline{L}$.  By Theorem 1A of \cite{hausel:2004}, the space $L^2\mathcal{H}^k$ of $L^2$ harmonic $k$-forms is
\begin{equation}
L ^2 \mathcal{H} ^k (K ) =
\begin{cases} 
H ^k (K , \partial K )\quad &\text{if $k < n + 1$},\\
\operatorname{Im} (H ^k (K , \partial K ) \rightarrow H ^k (K) ) \quad &\text{if $k = n + 1$},\\
H ^k (K )\quad &\text{if $k > n + 1$},
\end{cases} 
\end{equation} 
where $H ^k (K , \partial K ) \rightarrow H ^k (K) $ is the inclusion map.
We have
\begin{equation}
H ^k (K) \simeq H ^k (\mathbb{C}  P ^{ n } ) = \begin{cases} 
\mathbb{R}  \quad &\text{if $k =0, 2,4,  \ldots , 2n $},\\
0 &\quad \text{otherwise},
\end{cases} 
\end{equation} 
and, using Poincar\'e duality,
\begin{equation}
\label{compactlysupportedcohoml} 
H ^k (K , \partial K ) \simeq H _c  ^k (K) \simeq H ^{2n + 2-k} (K  ) = \begin{cases} 
\mathbb{R}  \quad &\text{if $k =2, 4, \ldots , 2n + 2$},\\
0 &\quad \text{otherwise}.
\end{cases} 
\end{equation} 
Finally the inclusion
\begin{equation}
\iota : H ^{n + 1} _c  (K ) \rightarrow H ^{n + 1} (K).
\end{equation} 
is the zero map for $n$ even, as  $H ^{n + 1} _c  (K ) =0$, and an isomorphism for $n$ odd.
In conclusion
\begin{equation}
L ^2 \mathcal{H} ^k (K ) = \begin{cases}
\mathbb{R}  \quad & \text{if } k =2, 4, \ldots , 2n,\\
0 &\text{otherwise} .
\end{cases} 
\end{equation} 
In particular, $L^2\mathcal{H}^2(K)$ is 1-dimensional.

We now show that we can write a generator $\tilde\omega  $ of $L ^2 \mathcal{H} ^2  (K) $ in the form $ \tilde\omega   = 2i\mathrm{d}  \beta  $ for $\beta  \in \Lambda ^{ 0,1 } (K) $ a Dirac zero mode. This extends the phenomenon observed in the case of the Calabi metric for $\mathbb{CP}^1$, where we had $\tilde{\omega} = 2i \mathrm{d} \beta$ with $\beta$ given by (\ref{asdfpot}). Furthermore, the explicit form (\ref{zm-sols}) of the zero modes of the twisted Dirac operator shows that $\beta$ is the $N=0$ zero mode of the untwisted Dirac operator, something which also holds for general $n$.

By the local $\partial \overline{\partial }$ lemma we can write $\tilde\omega   = 2i\partial\overline{\partial }  \phi $ for some locally defined real function $\phi$.  Set
\begin{equation}
\beta =\overline{\partial }\phi .
\end{equation} 
Clearly $\overline{\partial }\beta  =0 $ so $\beta $ is a Dirac zero mode provided that
\begin{equation}
\begin{split} 
0=D \beta  =- \overline{\partial }^\ast \beta =-\overline{\partial }^\ast \overline{\partial } \phi 
=g ^{ \mu \bar\nu } \partial _\mu \partial _{ \bar\nu } \phi = \overline{\partial }\partial \phi =\triangle \phi  .
\end{split} 
\end{equation} 
%If $D\beta =0 $  then  $\mathrm{d} \beta  $ is harmonic for
%\begin{equation}
%\mathrm{d} ^\ast \mathrm{d}\beta =  (\partial ^\ast + \overline{\partial }^\ast ) \partial \beta 
%= - \partial  (\partial ^\ast + \overline{\partial }^\ast ) \beta =0
%\end{equation} 
%since $\partial ^\ast |_{ \Lambda ^{ 0,k } }=0 $. 
Take $\phi$ to be a real function of $s$ and set $\chi =\mathrm{d} \phi / \mathrm{d} s $, so that
\begin{equation}
\beta = \chi (s)  \overline{\partial  } s .
\end{equation} 
We have $\partial \beta = \chi ^\prime \partial s \wedge \overline{\partial }s + \chi \partial \overline{\partial }s $ so contracting with the inverse metric we get
%\begin{equation}
%\partial \mathcal{A} = \chi ^\prime \partial s \wedge \overline{\partial }s + \chi \partial \overline{\partial }s,
%\end{equation} 
\begin{equation}
 - \overline{\partial }^\ast \beta 
 = g ^{\nu\bar \mu   } \left( \chi ^\prime \bar{z} _\mu z _\nu 
  + \chi \delta _{ \nu \mu } \right).  
\end{equation} 
%\begin{equation}
%\begin{split} 
% - \overline{\partial }^\ast A &
% = g ^{\nu\bar \mu   } \left( \chi ^\prime \bar{Z} _\mu Z _\nu 
%  + \chi \delta _{ \nu \mu } \right)  \\
%P ^- (\partial S \wedge \overline{\partial }S ) &
%= \sum _{ \mu, \nu }\bar{Z} ^\mu Z ^\nu P ^- (  \mathrm{d} Z ^\mu \wedge \mathrm{d} \bar Z ^\nu )
%= \sum _{ \mu, \nu }\bar{Z} ^\mu Z ^\nu g ^{\nu\bar \mu   },\\
%P ^- (\partial \overline{\partial }S )&
%= P ^- ( \sum _{  \mu  }\mathrm{d} Z ^\mu \wedge  \mathrm{d} \bar Z ^\mu ) = \sum _\mu  g ^{\bar \mu \mu } =\operatorname{Tr} (g ^{-1} ).
%\end{split}
%\end{equation} 
Using (\ref{ginverse})  we calculate
\begin{equation} 
\label{compsome} 
\begin{split} 
\operatorname{Tr} (g ^{-1} )&
= \frac{1}{F } \left( n + 1 + \frac{\kappa  }{s ^{n + 1 } } \right) 
= \frac{n}{F }   + F ^{ n  }  ,\qquad 
\bar{z} _\mu z _\nu g ^{ \nu  \bar\mu } 
 =\frac{1}{F} \left[ s + \frac{\kappa  }{s ^{ n   }  }  \right] = s F ^n  .
\end{split} 
\end{equation} 
Therefore, $\beta $ is a zero mode if
\begin{equation}
\chi ^\prime \left[ s + \frac{\kappa  }{s ^{ n }  }  \right]  + \chi \left( n + 1 + \frac{\kappa }{s ^{n + 1 } } \right)  =0,
\end{equation} 
which integrates to
\begin{equation}
\chi =\frac{1 }{s ^{n + 1 } \left( 1 + \frac{\kappa }{s ^{n + 1 } } \right) ^{ \frac{n}{n + 1 }} }
 =\frac{1}{s ^{n + 1 }  F^{ n} }.
\end{equation} 
Thus
\begin{align}\label{beta-general}
\beta  &=\chi  (s) z  \mathrm{d} \bar z
=  \frac{z  \mathrm{d} \bar z}{s ^{n + 1 }  F^{ n} },\\
\label{tilomega-general}
\tilde\omega &= 2i \mathrm{d} \left(\frac{z  \mathrm{d} \bar z}{s ^{n + 1 }  F^{ n} } \right).
\end{align} 
Note that for $n =1 $ we recover (\ref{asdfpot}).

The 2-form $\tilde \omega=2i\mathrm{d} \beta   $ has components
\begin{equation}
\tilde \omega _{ \mu \bar\mu } = - 2i\frac{n }{s ^{ n + 1 } F ^{ 2n + 1 }}, \qquad 
\tilde \omega _{ \mu \bar\nu } = -2i\left( 1 + \frac{n }{F ^{ n + 1 } } \right)  \frac{z _\mu \bar z _\nu  }{s ^{ n + 2 } F ^{ 2n + 1 }},\ \mu \neq \nu .
\end{equation} 
One can see that while $\beta $ is not defined for $s =0 $, $\tilde \omega $ is a globally defined 2-form.   To check the $L ^2 $ condition we need to integrate $|\tilde \omega |^2 $  with respect to the volume element $r ^{ 2n + 1 }\mathrm{d} r  $ where $r =\sqrt{ s }$. For large $r$ we have   $F \sim 1 $, $g ^{  \rho \bar\nu } \sim\delta ^{ \rho \nu } $, $|\tilde \omega |^2 \sim r ^{4 (n + 1 )}$ hence $|\tilde \omega  |^2 $  is integrable. On the other hand $\beta $ is not $L ^2 $. In fact since
\begin{equation}
\label{hodgedualA} 
\partial _{z_\mu }^\flat = F ^{  \prime } \bar{z} _\mu z _\nu \mathrm{d} z ^{ \bar\nu }+F \mathrm{d} z ^{ \bar\nu } \quad   \Rightarrow \quad 
z^\mu \partial _{z_\mu} ^\flat = \frac{{z} _\nu \mathrm{d} \bar z ^\nu }{F^{ n}},
\end{equation} 
the metric dual $\beta  ^\sharp $ of $\beta $   is 
\begin{equation}
\beta  ^\sharp =\beta  ^\nu \partial _{ z _ \nu }, \qquad 
{\beta }   ^\nu =  g ^{  \bar\mu \nu }{\beta } _{ \bar \mu } 
=\chi\,  z ^\nu F^{ n }
= \frac{ z ^\nu }{ s ^{n + 1} }.
\end{equation} 
Thus, the squared norm of $\beta$ is
\begin{equation}
|\beta  |^2 ={\beta } ^\nu \bar \beta  _{\nu} 
= \frac{1}{s ^{ 2n + 1 } F ^n } ,
\end{equation} 
which  is not $L ^2 $ due to the logarithmic divergence near $r =0 $.  

Finally, we want to normalise $\tilde \omega  $ so that it is the curvature of a connection $\mathcal{A}$. To do so,  we need to impose the quantisation condition
\begin{equation}
\frac{i}{2 \pi } \int _{ \Sigma } \mathrm{d} \mathcal{A} =\NN \in \mathbb{Z}  ,
\end{equation} 
where $\Sigma$ is any generator of $H ^2 (K , \mathbb{Z}  )=\mathbb{Z}  $. We can take $\Sigma$ to be the $\mathbb{C}  P ^1 $ obtained setting $ \zeta =0 =w _2 =\cdots =w _n $. To compute the flux of $\tilde \omega  $ over $\Sigma$ we switch to  $(w _i , \zeta )$ coordinates. Using (\ref{relzbardzeta}) we find
\begin{equation}
\begin{split}
\tilde \omega =\frac{2i}{ n + 1} \left( \frac{\mathrm{d} \log \zeta + \alpha  }{ (\kappa + |\zeta |^2  )^{ \frac{n}{ n + 1}}}  \right) ,
\end{split} 
\end{equation} 
for $\alpha$ given by (\ref{alphaform}). It follows that
\begin{equation}
\tilde \omega |_\Sigma = \frac{2i \mathrm{d}  \alpha |_{ \Sigma }}{(n + 1 ) \kappa ^{ \frac{n}{n + 1 }}}.
\end{equation} 
Since $ \tfrac{i}{2 \pi }\tfrac{\mathrm{d} \alpha  |_ {\Sigma }}{(n + 1 )}$ has unit flux over $\Sigma$,  the required normalisation is 
\begin{equation}
\label{genconnect} 
\mathcal{A} = \NN (A - \bar A ), \qquad  A =\kappa ^{ \frac{n}{n + 1} } \frac{z \mathrm{d} \bar z }{2s ^{ n + 1 }F ^n },
\end{equation} 
which for $n =1 $ gives back (\ref{twistingconn}).

The analogue of the $n =1 $ Killing vector field $ X _3 $ generating translation along the circles $ |\zeta  | = \mathrm{const} $ of the fibres is
\begin{equation}
\xi = \frac{i}{2}   ( z ^\nu \partial _\nu - \bar z ^\nu \partial _{ \bar\nu }).
\end{equation} 
%The hypersurfaces $ S =c>0 $ in $\mathcal{O} (-n )$ are the circles of the Hopf fibration 
%\begin{equation}
% S ^1 \hookrightarrow \mathscr{O} (-n-1 ) |_{S =\mathrm{const} }  \rightarrow \mathbb{C}  P ^{ n }
% \end{equation} 
%and we recognise (\ref{xikvf}) for the KVF generating the isometry given by translation along the fibres.
Using (\ref{hodgedualA}) we see that
\begin{equation}
 \xi  ^\flat = \frac{i}{2} \left( \frac{z \mathrm{d} \bar z - \bar z \mathrm{d} z }{F ^n }\right) \quad \Rightarrow \quad 
 2 \mathrm{d} \xi ^\flat 
= 2i  \mathrm{d} \left( \frac{ z \mathrm{d} \bar z }{ F ^n }\right) .
\end{equation} 
Since $F$ satisfies the identity
\begin{equation}
F - \frac{\kappa }{F ^n s ^{ n + 1 }} = \frac{1}{F ^n },
\end{equation} 
it follows that
\begin{equation}\label{X-dual-omega}
\omega - \kappa \tilde \omega 
=2 i \mathrm{d}\left( \frac{  z \mathrm{d} \bar z }{F ^n } \right) =2 \mathrm{d} \xi ^\flat .
\end{equation} 

%In fact the total space is actually a $\mathbb{Z}  _n $ quotient of that, where $\zeta \in  \mathbb{Z}  _n $ acts diagonally on the homogeneous coordinates $z ^i $, $\zeta \cdot (z ^1 , \ldots , z ^n )= (\zeta z ^1 , \ldots , \zeta z ^n )$. Thus the $\mathbb{Z}  _n $ quotient leaves the $\mathbb{C}  P ^{ n-1 }$ base invariant, and on the complex line bundle plane leaves the radius invariant and shortens the fibre angle range from $4 \pi $ to $4 \pi /n $. 
%Since the manifold is Ricci-flat, $\mathrm{d} \xi ^\flat $ is harmonic for any KVF $ \xi $, which we already knew,  but clearly $\xi$ has a special role being the isometry group of the fibre. It is also special since it is $\Xi $ is the Green's function for the Laplacian.

\subsection{Zero modes of the twisted Dirac operator on $\mathcal{O}(-n-1)$}
As already mentioned, by the same argument used for  EH,  the canonical spin-$c$ structure on $K=\mathcal{O}(-n-1)$  admits no non-trivial $L ^2 $ harmonic spinor. In this section we consider the problem of normalisable zero modes with respect to the Dirac operator twisted by the  connection (\ref{genconnect}). 
We will limit ourselves to the study of zero modes in $\Lambda ^{ 0,1 } (K)  $.

Take the ansatz
\begin{equation}\label{zm-general}
\sigma = P (z _i ) h (s) \overline{\partial }s,
\end{equation} 
where $P$ is a function of $z_1, \ldots, z_{n+1}$ only,  so that  the equation $\overline{\partial }\sigma =0 $ is automatically satisfied. The remaining equation  is the projection onto $\Lambda ^{ 0 } (K) $ of 
\begin{equation}
\partial   \sigma = -\NN \bar{A }\wedge \sigma .
\end{equation} 
One has
\begin{equation}
\partial \sigma =h \partial _i P\, \mathrm{d} z ^i\wedge  \overline{\partial }s  + P h ^\prime  \partial s \wedge \overline{\partial }s + Ph \partial \overline{\partial }s.
\end{equation} 
The projections of  $\partial \overline{ \partial } s$, $ \partial s \wedge   \overline{\partial } s $ onto $\Lambda^0(K)$ are given by (\ref{compsome}), and
\begin{equation}
h z _\mu \partial _\nu P     g ^{  \mu \bar\nu }
= h z  _\mu \partial _\nu f \left(  \frac{\delta ^{ \nu\mu }}{F } + \frac{\kappa  }{s ^{ n  + 2 } F  }\bar z ^\mu  z ^\nu  \right)
=h F ^n z _\mu \partial _\mu P. 
\end{equation} 
Now take $P$ to be a homogeneous polynomial in $ (z _1 , \ldots , z _{ n + 1 } )$ of degree $\dd $, so that $ z _\mu \partial _\mu P = \dd P$. Then
\begin{equation}
- \overline{\partial } ^\ast \sigma 
= \frac{ P}{F} (F ^{ n + 1 }  ( \dd h  + h ^\prime s + h ) + n h).
\end{equation} 

We also need to compute the $\Lambda^0(K)$ projection of $\mathcal{A} \wedge \sigma = \bar A\wedge \sigma$.
Since 
\begin{equation}
\bar{A } = \kappa ^{ \frac{n}{n + 1} } \frac{ \partial s }{2s ^{ n + 1 }F ^n } ,
\end{equation}  we obtain
\begin{equation}
 - \bar{A }_\mu  \sigma _\nu g ^{ \mu \bar\nu }
 = - \kappa ^{ \frac{n}{n + 1} } \frac{ P h z _\nu \bar z _\mu g ^{ \nu  \bar\mu }  }{2s ^{ n + 1 }F ^n }
 =- \kappa ^{ \frac{n}{n + 1} } \frac{ Ph }{2s ^{ n  } }.
\end{equation} 
Therefore, the twisted Dirac equation becomes the following ODE for $h$,
\begin{equation} 
\label{twisdirac} 
F ^{ n + 1 }  ( (\dd + 1 ) h  + h ^\prime s  ) + n h
 = \kappa ^{ \frac{n}{n + 1} }\NN  \frac{hF  }{2s ^{ n  } }.
\end{equation} 
Equivalently
\begin{equation} 
\label{odeforh}
  (\log h )^\prime     
 = \frac{\NN \kappa ^{ \frac{n}{n + 1} }}{2  F ^{ n  } s ^{ n + 1  } } -   \frac{ \dd + 1}{s} - \frac{n}{s F ^{ n + 1 } },
\end{equation} 
which has solution
\begin{equation}
h = h _0f, \qquad  h _0  = \frac{1}{F ^n s ^{ \dd + n + 1 }}.
\end{equation} 
Here $h_0$ solves (\ref{odeforh}) for $\NN=0$, and $f$ satisfies
\begin{equation}
(\log f )^\prime 
= \frac{\NN \kappa ^{ \frac{n}{n + 1} }}{2  F ^{ n  } s ^{ n + 1  } }.
\end{equation} 

To discuss normalisability we only need the behaviour of $f$  for small and large $s$, which is
\begin{equation}
f \sim \begin{cases}
s ^{ \NN /2} & \text{for $s\ll 1  $},\\
\exp \left(  - \frac{\NN \kappa ^{ \frac{n}{n + 1 }}}{2n} \frac{1 }{s ^n } \right) & \text{for $s\gg1  $}.
\end{cases} 
\end{equation} 
Since $|P |^2 =s ^{ \dd} $, $|\overline{\partial }s|^2  = F ^n s $, the squared norm of $\sigma$ is
\begin{equation}
|\sigma |^2 = |P |^2  h ^2  |\overline{\partial } s |^2 
%= s ^{ \dd} \frac{1}{F ^{ 2n } s ^{ 2 \dd  + 2n + 2 }} s F ^n f ^2 
= \frac{  f ^2 }{F ^{ n } s ^{  \dd  + 2n + 1 }} ,
\end{equation} 
which is to be integrated with respect to the volume element $r ^{ 2n + 1 }$ for $r =\sqrt{ s }$.
For large $r$  we have $F \sim 1 $ so $|\sigma | ^2 r ^{2n + 1 } \mathrm{d} r \sim r ^{ - (2 \dd + 2n + 1)} \mathrm{d} r $  which always gives a  finite contribution. For small $r$ we have $f \sim r^{ \NN } $, $F \sim \kappa ^{ \frac{1}{n + 1 }}r ^{-2} $, hence   $|\sigma | ^2 r ^{2n + 1 } \mathrm{d} r \sim r ^{ 2\NN - 2\dd -1 }\mathrm{d} r $ and square-integrability requires 
\begin{equation}
\label{genl2cond} 
\NN > \dd.
\end{equation} 
Note the similarity with the condition (\ref{ehl2cond}) obtained for $n =1 $, where we took $\dd=2 N $. Note also how the (untwisted, not $L^2$)  spinor $h_0 \,\overline \partial s$ obtained for $\NN=0$, $\dd=0$ is, up to scale, equal to the connection (\ref{genconnect}) with harmonic $L^2$ curvature, again in complete analogy with the $n=1$ case.

\section{Discussion}

The main result of this work is the explicit description of the $L^2$ zero modes of the (twisted) Dirac operator on both the Eguchi-Hanson metric, see (\ref{zm-sols}),  and its higher dimensional generalisation given by the Calabi metric on $\mathcal{O}(-n-1)$. 
As expected, for $n=1$ the EH zero modes organise themselves into multiplets of the EH isometry group $SU(2)$. The dimension of the space of zero modes is controlled by the integer $\NN$, see (\ref{twistingconn}), that controls the twist. The Dirac operator index analysis in \cite{franchetti:2018} confirms that the zero modes obtained are all the zero modes, and so for the EH space the problem of finding the harmonic spinors is completely solved.

Far from arbitrary, the $U(1)$ connection used for the twist is preferred geometrically for multiple reasons. First, its curvature is the unique (up to scale) harmonic $L^2$ 2-form on EH. Second, it is a connection on the total space of the $S^1$ bundle over EH arising as a level set in the process of the hyperk\"ahler reduction from $\mathbb{H}^2 $. In passing, this suggests that the twisted Dirac operator may also be understood as an appropriate dimensional reduction of the untwisted Dirac operator on $\mathbb{H}^2$. It would be very interesting to see whether this is the case, and whether  the harmonic spinors (\ref{zm-sols}) can also be understood as arising from the dimensional reduction of some very simple spinors on $\mathbb{C}^4$. We leave this to further work.
Third, the $U(1)$ connection with $L^2$ curvature agrees with the lowest lying (untwisted, not $L^2$) Dirac zero mode.

The EH metric is the $n=1$ case of Calabi's family of  Ricci-flat K\"ahler metrics on $\mathcal{O}(-n-1)$, $n\in\mathbb{N}$. As we have shown, many of the EH results generalise to higher $n$. First, the  description (\ref{calabimetriccpn}) of the metric in terms of "symmetrical" coordinates   completely parallels the case of EH. Second, for all values of $n$ there is a unique $L^2$ harmonic 2-form $\tilde\omega$, see (\ref{tilomega-general}), which also arises as the curvature of the lowest lying zero mode of the untwisted Dirac operator, and whose connection can be used to twist the Dirac operator. Moreover, $\tilde\omega$ differs from the K\"ahler form $\omega$ by a constant multiple of $\mathrm{d}\mathrm\xi^\flat$, see (\ref{X-dual-omega}), for $\xi$ a Killing vector field generating the isometric $U(1)$ action on the fibre. Interestingly, the forms $\omega$, $\tilde\omega$, $\mathrm{d}\xi^\flat $ are all harmonic although only $\tilde\omega$ is $L^2$.  Third,  the EH zero modes (\ref{zm-sols}) of the Dirac operator twisted by this preferred $U(1)$ connection have analogues (\ref{zm-general}) in the general case. These general zero modes  again fall into irreducible representations of the $U(n+1)$ isometry group, and the number of allowed zero modes is controlled by the integer $\NN$ that determines the curvature flux. 

The cases $n=1$ and $n>1$ also present an important difference: the EH metric is hyperk\"ahler while for $n>1$  Calabi's metric is Calabi-Yau but not hyperk\"ahler.  For this reason, for $n>1$ there can be no analogue of the hyperk\"ahler quotient derivation of the EH metric reviewed in section \ref{subsec:hyperkquot}, but one may still wonder if the metric could be obtained as the K\"ahler reduction of some simpler higher-dimensional metric.

In the general case, we cannot claim that the Dirac zero modes that we have found exhaust all $L^2$ harmonic spinors. For EH it is possible to come to this conclusion thanks to the fact that, since the curvature of the twisting connection is self-dual, there are no zero modes in $W^+$ and the index of the Dirac operator is equal to  the number of $L^2$ zero modes in $W^-$.  This self-duality argument does not extend to higher $n$, and we have no alternative argument  showing that there are no $L^2$ zero modes in $W^+$.  In fact, we have only analysed  zero modes belonging to $\Lambda^{0,1}$ and we also do not know if there are  additional zero modes in the spaces $\Lambda^{0,2k+1}$, $2k\leq n$. A complete answer to these questions is left for future work, along with the interesting problem of studying zero modes of the Dirac operator twisted by $L^2$ harmonic forms of degree other than two.

Another outcome of this work is the development of a set of rules for computing the spin-$c$ Dirac operator on a K\"ahler manifold. Indeed, we have shown that the computation of the action of this Dirac operator on spinors is no more complicated than the computation of the exterior derivative of differential forms. The only additional operation needed is the application of the metric contraction to the result of the exterior derivative, to map the latter into the space $\Lambda^{0,\bullet}$ where  spinors live. This gives extremely simple computational rules  and makes the formalism of spinors as differential forms extremely convenient for explicit calculations involving the Dirac operator. We hope this work will lead to a better familiarity of the community with the very efficient computational tool that spin-$c$ spinors provide.

We close with some remarks on what motivated us to embark on the present investigation. The Calabi construction can be applied to an arbitrary K\"ahler manifold $M$ with non-zero scalar curvature, and gives a Ricci-flat K\"ahler metric on the total space of the canonical bundle of $M$. Applying this construction to  $\mathbb{CP}^1\times\mathbb{CP}^2$ is particularly interesting because the resulting Ricci-flat metric has the Standard Model gauge group as its isometry group. The resulting  space is Calabi-Yau with  holonomy group $SU(4)$. Its metric is asymptotically conical, with the metric on the base of the cone being the nearly parallel $G_2$ metric $M(3,2)$ discussed in particular in \cite{FKMS}. We are interested in determining the $L^2$ harmonic spinors of (the appropriately twisted) Dirac operator on this space. The present work can be considered as setting the stage for this more involved computation. 

\subsection*{Acknowledgements} The authors are grateful to Bernd Schroers for introducing them to each other, and for participating in the early stages of this collaboration. 
GF thanks the Simons Foundation for its support under the Simons Collaboration on Special Holonomy in Geometry, Analysis and Physics [grant number 488631].

\appendix
\section{The geometry of $\mathbb{C}  ^2 $}
\label{c2stuff} 
In this appendix we want to relate the descriptions of  Euclidean $\mathbb{R}  ^4  =\mathbb{C}  ^2 $ as a cohomogeneity-one space with respect to the action of $SU (2) $, as a complex manifold and as a complex line bundle. To that end it is convenient to first review the geometry of $SU (2) $.

\subsection{The geometry of $SU (2) $}
As it is well known, $S ^3 \simeq SU (2) $. A possible parametrisation of left-invariant 1-forms on $SU (2) $ is 
\begin{equation}
\label{livforms} 
\begin{split} 
\eta _1 &= +  \sin \psi  \, \mathrm{d} \theta - \cos \psi \sin \theta \, \mathrm{d} \phi ,\\
\eta _2 &= -  \cos \psi \, \mathrm{d} \theta - \sin \psi \sin \theta \, \mathrm{d} \phi  ,\\
\eta _3 &= \mathrm{d} \psi + \cos \theta \, \mathrm{d} \phi,
\end{split} 
\end{equation} 
where $\theta \in [0, \pi ] $, $\phi \in [0, 2 \pi )$, $\psi \in [0, 4 \pi )$. 
Taking the same range for $\theta, \phi $ with $\psi \in [0, 2 \pi )$ gives instead a parametrisation of $SU (2) / \mathbb{Z}  _2 =SO (3)  $.
Note that the forms (\ref{livforms}) satisfy
\begin{equation}
\mathrm{d} \eta _i = +\frac{1}{2}  \epsilon _{ ijk } \eta _j \wedge \eta _k .
\end{equation} 

The left-invariant vector fields on $SU (2) $  satisfying $\eta _i (X _j )=\delta _{ ij }$ are
\begin{equation}
\label{livvecfields} 
\begin{split}
X _1 &=+  \sin \psi\,  \frac{ \partial }{\partial \theta}  + \frac{\cos \psi }{\sin \theta } \left( \cos \theta \, \frac{\partial }{\partial \psi } - \frac{\partial }{\partial \phi } \right) ,\\
X _2 &= - \cos \psi\,  \frac{\partial }{\partial \theta} +  \frac{\sin  \psi }{\sin \theta } \left( \cos \theta \, \frac{\partial }{\partial \psi } - \frac{\partial }{\partial \phi } \right) ,\\
X _3 &=+ \frac{ \partial }{ \partial\psi}.
\end{split}
\end{equation} 
They satisfy the $\mathfrak{su }(2) $ Lie algebra relation
\begin{equation}
[ X _i , X _j ] = - \epsilon _{ ijk } X _k .
\end{equation} 
In fact it can be checked that $X _i $ is the left-invariant vector field associated to $ \tfrac{i}{2} \sigma _i $, for $(\sigma _i )$ the Pauli matrices.

We can also introduce right-invariant 1-forms $(\zeta _i )$ and vector fields $( Y _i )$, 
\begin{equation}
\label{rivforms} 
\begin{split} 
\zeta  _1 &= - \cos \phi \sin \theta  \, \mathrm{d} \psi  +  \sin \phi \, \mathrm{d} \theta ,\\
\zeta  _2 &= + \sin \phi \sin \theta \, \mathrm{d} \psi  + \cos \phi  \, \mathrm{d} \theta   ,\\
\zeta  _3 &= \mathrm{d} \phi + \cos \theta \, \mathrm{d} \psi,
\end{split} 
\end{equation} 
\begin{equation}
\label{rivvecfields} 
\begin{split} 
Y _1 &= \sin \phi \frac{\partial }{\partial \theta } + \frac{\cos \phi }{\sin \theta } \left( \cos \theta \frac{\partial }{\partial \phi } - \frac{\partial }{\partial \psi } \right),\\ 
Y _2 &= \cos  \phi \frac{\partial }{\partial \theta } - \frac{\sin  \phi }{\sin \theta } \left( \cos \theta \frac{\partial }{\partial \phi } - \frac{\partial }{\partial \psi } \right),\\ 
Y _3 &=\frac{\partial }{\partial \phi },
\end{split} 
\end{equation} 
satisfying $\zeta _i (Y _j ) = \delta _{ i j }$,
\begin{equation}
[ Y _i , Y _j ] = + \epsilon _{ ijk }Y _k , \quad \mathrm{d} \zeta _i =-  \frac{1}{2} \epsilon _{ ijk }\zeta _k , \quad [X _i , Y _j ]=0.
\end{equation}

The metric 
\begin{equation}
\label{s3metric} 
g _{ S ^3 }  =  \frac{1}{4} (\eta _1 ^2 + \eta _2 ^2 + \eta _3 ^2 )=  \frac{1}{4} (\zeta  _1 ^2 + \zeta  _2 ^2 + \zeta  _3 ^2 )
\end{equation} 
is the round metric on the 3-sphere of unit radius, or equivalently the bi-invariant metric on $SU (2) $. In terms of the latter description it is clear that the metric is invariant under both a left and right $SU (2) $ action. The left action is generated by the \emph{right}-invariant vector fields  (\ref{rivvecfields}). Such vector fields satisfy $ L _{ Y _i } \eta  _j =0 $. The right $SU (2) $ action is generated by the \emph{left}-invariant vector fields (\ref{livvecfields}). We have
\begin{equation}
\begin{split} 
L _{ X _i } \zeta _j &= L _{ Y _i } \eta _j =0,\\
L _{ X _i } \eta _j &=- \epsilon _{ ijk } \eta _k, \quad L _{ Y _i } \zeta _j = +  \epsilon _{ ijk } \zeta _k .
\end{split} 
\end{equation} 

\subsection{The geometry of $\mathbb{C}  ^2 $}
The orbit structure of $\mathbb{C}  ^2 \simeq \mathbb{R}  ^4  $ with respect to the $SU (2) $ action is obtained by writing $ \mathbb{R}  ^4=  \{ 0 \}\cup \big( (0, \infty ) \times S ^3 \big)   $. Introducing a radial coordinate $r \in [0, \infty )$ transverse to the $SU (2) $ orbits we can write the Euclidean metric on  $\mathbb{C}  ^2 $  as
\begin{equation}
\label{su2metricc2} 
\gc= \mathrm{d} r ^2 + \frac{r ^2 }{4} (\eta _1 ^2 + \eta _2 ^2 + \eta _3 ^2 ).
\end{equation} 
The quantities  $(\eta _i ) $, $(\zeta _i )$, $(X _i )$, $(Y _i )$  extend to  well-defined 1-forms and vector fields on 
$\mathbb{R}  ^4 \setminus \{ 0 \} $.

As a real manifold $\mathbb{R}  ^4 $ has isometry group $ O (4) $. We can identify a point $(z _1 , z _2 ) \in \mathbb{C}  ^2 \setminus \{ 0 \} $ with a pair $(r, x )$ where  $x$ is  the $SU (2) $ element given by
\begin{equation} \frac{1}{|z _1 |^2 + |z _2 |^2 }
\label{c2assu2} 
\begin{pmatrix}
 z _1  & - \bar z _2  \\
z _2  &  \bar z _1 
\end{pmatrix} 
\end{equation} 
and $ r = \sqrt{ |z _1 |^2 + |z _2 |^2 }$. Left and right matrix multiplication of (\ref{c2assu2}) by an $SU (2) $ element give raise to isometries, with the left action generated by the right-invariant vector fields $(Y _i )$ and the right action by the left-invariant vector fields $(X _i )$.  Explicitly is $h \in SU (2) $ is given by
\begin{equation}
h = \begin{pmatrix}
 a & - \bar b \\
b & \bar a
\end{pmatrix} 
\end{equation}  
the left and right $SU (2) $ action are given by (\ref{su2la}), (\ref{su2ra}) respectively,
\begin{align} 
\label{su2la} 
z _1 &\mapsto a z _1 - \bar b z _2 , \quad  z _2 \mapsto b z _1 + \bar a z _2 ,\\
\label{su2ra} 
z _1 &\mapsto a z _1 - b \bar z _2 , \quad  z _2 \mapsto a z _2 + b \bar  z _1 .
\end{align} 
Since $(SU (2) \times SU (2)) / \mathbb{Z}  _2  =SO (4) $, the left and right $SU (2) $ action give  the (connected component of the) full isometry group  of $\mathbb{R}  ^4 $.  The $\mathbb{Z}  _2 $ quotient corresponds to the fact that for $a =\bar a $, $b =0 $, which in $ SU (2) $ implies $a =\pm 1 $, left and right action result in the same transformation.

The map (\ref{su2ra}) is not complex linear. In fact $\mathbb{C}  ^2 $ as a complex manifold has metric
\begin{equation}
\label{complmetricc2} 
\gc= |\mathrm{d} z _1 | ^2 + |\mathrm{d} z _2 |^2 ,
\end{equation} 
 complex structure 
\begin{equation} 
J (\partial / \partial z _i )= + i \partial / \partial z _i  , \quad J (\mathrm{d} z _i )=- i \mathrm{d} z _i ,
\end{equation} 
and   the smaller isometry group $U (2) \subset SO (4) $. 
%Besides the left acting $SU (2) $, there is a remaining diagonal $U (1) $ action 
%the remaining isometries are generated by $( X _i )$ and correspond to the right action of $SU (2) $ on $\mathbb{C}  ^2 $.
There is a group isomorphism $U (2) = (SU (2) \times U (1) )/ \mathbb{Z}  _2  $, $(h , u ) \mapsto h u $. The (left) $SU (2) $ action is (\ref{su2la})  and the  $U (1) $ action is diagonal,
\begin{equation}
\label{u1action} 
 (z _1 , z _2 ) \cdot \mathrm{e} ^{ it }=   (z _1 , z _2 ) \mathrm{e} ^{ it }.
\end{equation} 
Looking at (\ref{su2ra}) we see that (\ref{u1action}) is the $U (1) $ action obtained restricting the  right  $SU (2) $ action to the $U (1) $ subgroup obtained by  setting $b =0 $.

%The $U (1) $ action is Hamiltonian with Hamiltonian function
%\begin{equation}
%\mu  =  \frac{1}{2} ( |z _1 |^2 +  |z _2 |^2 ).
%\end{equation} 

We now define an orthonormal frame $(X _1 , X _2 , X _3 , X _4 )$  consisting of real vector fields adapted to the action (\ref{u1action}) and the complex structure $J$,
\begin{equation}
\label{adaptedframe} 
\begin{split} 
X_1 & =- \frac{i}{2}  \left(  \bar{z} _2 \frac{\partial }{\partial z _1 } - \bar{z} _1 \frac{\partial }{\partial z _2 } - z _2 \frac{\partial }{\partial \bar z _1 }+  z _1\frac{\partial }{\partial \bar z _2} \right) ,\\
X _2 &=-  \frac{1}{2} \left(\bar{z} _1 \frac{\partial }{\partial z _2 } - \bar{z} _2 \frac{\partial }{\partial z _1 } 
+ z _1 \frac{\partial }{\partial \bar z _2 }- z _2 \frac{\partial }{\partial \bar z _1 }\right) , \\
X_3 &=+  \frac{ i}{2} \left(  z _1 \frac{\partial }{\partial z _1 } + z _2 \frac{\partial }{\partial z _2 } -\bar{z}  _1 \frac{\partial }{\partial \bar{z}  _1 }-\bar{z}  _2 \frac{\partial }{\partial \bar{z}  _2 }\right) ,\\
X _4 & =-\frac{1}{2}  \left(  z _1 \frac{\partial }{\partial z _1 } + z _2 \frac{\partial }{\partial z _2 } + \bar{z}  _1 \frac{\partial }{\partial \bar{z}  _1 } + \bar{z}  _2 \frac{\partial }{\partial \bar{z}  _2 }\right) .
\end{split} 
\end{equation} 
The frame $(X _i ) $ is adapted in the sense that $X _3 $ is the infinitesimal generator of the action (\ref{u1action}), $X _4 =J X _3 $, $X _2 =J X _1 $  and 
$(X _1 , J X _1 )$ is  the $\gc $-orthogonal complement of $(X_3, JX_3) $. More precisely,
\begin{equation}
\label{u1gen} 
X _3 -i J X _3  = i \left(  z _1 \frac{\partial }{\partial z _1 } + z _2 \frac{\partial }{\partial z _2 } \right) 
\end{equation} 
generates (\ref{u1action}) and $X _3 $ is the corresponding real holomorphic vector field.
 We have
\begin{equation}
|X _1 |^2 =|X _2 |^2= |X _3 |^2  =|X _4 |^2  = \frac{1}{4} (|z _1 |^2 + |z _2 |^2 ),
\end{equation} 
where $ |X_i |^2 =\gc(X_i,X_i) $.
 As we will show, the vector fields $(X _1 , X _2 , X _3 )$ in (\ref{adaptedframe}) are the left-invariant vector fields  (\ref{livvecfields}) expressed in terms of complex coordinates. 

Denoting by $ X ^\flat $ the $\gc $-metric dual of $X$ we have
\begin{equation} 
\label{adaptedcoframe} 
\begin{split}
\theta _1 &= X _1 ^\flat 
=- \frac{ i}{4} \left(  z _1 \mathrm{d} {z} _2 - z _2 \mathrm{d} {z} _1  -\bar{z}  _1  \mathrm{d} \bar z _2  + \bar{z}  _2 \mathrm{d} \bar z _1 \right) 
= - \frac{1}{2} \operatorname{Im} (z _1 \mathrm{d} z _2 - z _2 \mathrm{d} z _1 ),\\
\theta _2 &= X _2 ^\flat 
=- \frac{ 1}{4} \left(  z _1 \mathrm{d} {z} _2 - z _2 \mathrm{d} {z} _1  + \bar{z}  _1  \mathrm{d} \bar z _2 -\bar{z}  _2 \mathrm{d} \bar z _1 \right) 
= + \frac{1}{2} \operatorname{Re} (z _1 \mathrm{d} z _2 - z _2 \mathrm{d} z _1 ),\\
\theta_3  &=X_3 ^\flat  
=+ \frac{ i}{4} \left(  z _1 \mathrm{d} \bar{z} _1  + z _2 \mathrm{d} \bar{z} _2  -\bar{z}  _1  \mathrm{d} z _1 -\bar{z}  _2 \mathrm{d} z _2 \right)
= + \frac{1}{2} \operatorname{Im} (\bar{z} _1 \mathrm{d} z _1 + \bar{z} _2 \mathrm{d} _2 ) ,\\
\theta _4  &=X _4 ^\flat  
=- \frac{1}{4}  \left(  z _1 \mathrm{d} \bar{z} _1  + z _2 \mathrm{d} \bar{z} _2   + \bar{z}  _1  \mathrm{d} z _1  + \bar{z}  _2 \mathrm{d} z _2 \right) 
= -\frac{1}{2}  \operatorname{Re}  (\bar{z} _1 \mathrm{d} z _1 + \bar{z} _2 \mathrm{d} _2 ) .
\end{split} 
\end{equation} 
Thus
\begin{equation}
\label{gc2split} 
\begin{split} 
\gc&
= \sum _{ i =1 }^4 \frac{ \theta _i ^2 }{ |X _i |^2 }
= \frac{  |\bar{z} _1  \mathrm{d} z _1 + \bar{z} _2 \mathrm{d} z _2 |^2 }{{|z _1 |^2 + |z _2 |^2 }} + \frac{ |z _1 \mathrm{d} z _2-  z _2 \mathrm{d} z _1 |^2 }{|z _1 |^2 + |z _2 |^2 },
\end{split} 
\end{equation} 
where the first and second addend correspond to the metric restricted to $ \operatorname{Span} (X_3,JX_3 )$, $ \operatorname{Span} (X_1,JX_1 )$ respectively.

The metrics (\ref{complmetricc2}), (\ref{su2metricc2})   are related by the diffeomorphism
\begin{equation}
\label{c2coords} 
\begin{split} 
z _1 &=r\cos \left( \frac{\theta }{2} \right)  \mathrm{e} ^{ \frac{i}{2} ( \psi + \phi  ) },  \quad 
 z _2 =r\sin \left( \frac{\theta }{2} \right)  \mathrm{e} ^{ \frac{i}{2} ( \psi - \phi  ) }.
 \end{split} 
\end{equation} 
It can be checked the pulling back by (\ref{c2coords}) the vector fields $( X _1 , X _2 , X _3 )  $ in (\ref{adaptedframe})  map to (\ref{livvecfields}), and
\begin{equation}
X _4 = - \frac{r}{2}  \frac{\partial }{\partial r}, \qquad \theta _4 =- \frac{r}{2} \mathrm{d} r.
\end{equation} 
Note how $(X _3, X _4 )$ span the $(r, \psi )$ plane.
The forms $\eta _i $, $\theta  _i $ are the metric dual of  $X _i $ with respect to, respectively,  $g _{ S ^3} $ and $\gc$. For  vector fields $X,Y$ tangent to $S ^3 $ we have $g _{ S ^3 } (X  ,Y) =  \tfrac{4}{ |z _1 |^2 + |z _2 |^2  }\gc( X,Y )$, hence
\begin{equation}
\eta _i =\left(  \frac{4}{|z _1 |^2 + |z _2 |^2 } \right) \theta _i .
\end{equation}

The space  $\mathbb{C}  ^2 \setminus \{ 0 \}  $ can be viewed as a line bundle over $\mathbb{C}  P ^1 $. We introduce  complex coordinates $(w, \zeta )\in \mathbb{C}  ^2 $ on the base and fibre via
\begin{equation} 
\label{compltobundle} 
(z_1, z _2 ) = \frac{\zeta }{\sqrt{ 1 + |w |^2 } } (w , 1 ), \qquad 
(w, \zeta ) =\left( \frac{z _1 }{z _2 }, \sqrt{ |z _1 |^2 + |z _2 |^2 } \frac{z _2 }{|z _2 |}\right) .
\end{equation} 
Since
\begin{equation}
\begin{split} 
|z _1 |^2 + |z _2 |^2 &= |\zeta |^2 ,\\
z _1 \mathrm{d} z _2 - z _2 \mathrm{d} z _1 &= - \frac{\zeta ^2   \mathrm{d} w}{1 + |w| ^2 },\\
\bar{z} _1 \mathrm{d} z _1 + \bar{z}  _2 \mathrm{d} z _2 &
= \bar\zeta \mathrm{d}\zeta +   \frac{  |\zeta |^2}{2} \left( \frac{ \bar{w} \mathrm{d} w - w  \mathrm{d} \bar{w} }{1 + |w| ^2 } \right) 
= \bar\zeta (\mathrm{d}\zeta +   i a\zeta),
\end{split} 
\end{equation} 
where
\begin{equation} 
\label{gc2bundle} 
a = \frac{1}{2i} \frac{\bar w \mathrm{d} w - w  \mathrm{d} \bar w}{1+|w|^2}
= \frac{ \operatorname{Im} (\bar w \mathrm{d} w )}{1+|w|^2},
\end{equation} 
we  get
\be
\label{flat-metric-zt}
\gc= |\mathrm{d} \zeta  +   \im a \zeta  |^2 + \frac{ |\zeta |^2 |\mathrm{d} w|^2}{(1+|w|^2)^2}.
\ee

Thinking of $S ^3 $ as the Hopf fibration $ U (1) \hookrightarrow S ^3\rightarrow  \mathbb{C}  P ^1 $ we recognise $z$ as an inhomogeneous coordinate on the base $\mathbb{C}  P ^1 $ and $\zeta $, for $|\zeta |  = |R| $, as the angle parametrising the $U (1) $ fibres. If the value of $|\zeta |$ is allowed to vary in $(0, \infty )$ we  obtain a parametrisation of $ (0, \infty ) \times S ^3 $ viewed as  the line bundle $ \mathbb{C}  ^2 \setminus \{ 0 \}  \rightarrow \mathbb{C}  P ^1 $.
Since 
\begin{equation}
\zeta \frac{\partial }{\partial \zeta } = z _1 \frac{\partial }{\partial z _1 } + z _1 \frac{\partial }{\partial z _2 }
\end{equation} 
the vector field $X _3 $ translating along the fibres becomes
\begin{equation}
X _3 = \frac{i}{2} \left( \zeta \frac{\partial }{\partial \zeta } - \bar \zeta \frac{\partial }{\partial \bar\zeta } \right) .
\end{equation} 

Further parametrising
\be
\label{zetar} 
\zeta = r e^{i\chi},  \quad \Leftrightarrow \quad 
 r = |\zeta |=\sqrt{ |z _1 |^2 + |z _2 |^2  }, \ \mathrm{e} ^{ i \chi } =  \frac{z _2 }{|z _2 |},
\ee
we obtain
\be\label{R4-zr}
\gc=  \mathrm{d} r^2 + r^2  \left[ ( \mathrm{d} \chi +  a)^2 + \frac{| \mathrm{d} w|^2}{(1+|w|^2)^2}\right],
\ee
where the term in brackets is the round metric on $S^3$. 
Introducing the usual coordinates  $( \theta , \phi )$, $\psi\in [0, 4 \pi )$ parametrising the base and fibre of the Hopf fibration as in (\ref{c2coords})  we get
\begin{align}
w &=\frac{z _1 }{z _2 }=  \cot \frac{\theta }{2} \mathrm{e} ^{ i \phi },\\
\label{chiandphisu2} 
 \chi &=\frac{ \psi -\phi }{2} ,
\end{align} 
so that
\begin{equation}
\label{sphercoordsc4} 
\begin{split} 
 \frac{|\mathrm{d} w|^2}{(1+|w|^2)^2} &
 =\frac{1}{4} (\mathrm{d} \theta ^2 + \sin ^2 \theta\,  \mathrm{d} \phi ^2 ) 
 =\frac{1}{4} (\eta _1 ^2 + \eta _2 ^2 ), \\
  \mathrm{d} \chi + a &= \frac{\mathrm{d} \psi + \cos \theta \, \mathrm{d} \phi }{2} =\frac{\eta _3 }{2}, \quad  a =\left(  \frac{1 + \cos \theta }{2} \right) \mathrm{d} \phi ,
 \end{split} 
\end{equation} 
 recovering (\ref{su2metricc2}).

The line bundle $ \mathbb{C}  ^2 \setminus \{ 0 \}  \rightarrow \mathbb{C}  P ^1 $ is homotopy equivalent to the Hopf fibration. The latter has first Chern number one. In fact
\begin{equation}
\beta=[\bar z _1 \mathrm{d} z _1 + \bar z _2 \mathrm{d} z _2 ]|_{ S ^3 }
= i \operatorname{Im} (\bar z _1 \mathrm{d} z _1 + \bar z _2 \mathrm{d} z _2 ) |_{ S ^3 }
= i ( \mathrm{d} \chi + a ) 
=\frac{i}{2} \eta _3 
\end{equation}  
evaluates to 1 on the generator (\ref{u1gen}) of the $U (1) $ action and so defines a connection on the total space of the Hopf bundle. By Chern-Weyl the first Chern number of the Hopf fibration is 
\begin{equation}
\frac{i}{2 \pi } \int _{ \mathbb{C}  P ^1 }  \mathrm{d} \beta   =1.
\end{equation}

The splitting (\ref{gc2split}) can be understood as a special case of a K\"ahler reduction.
Let $ M$ be a K\"ahler manifold with metric $g$ and K\"ahler  form $\omega$. Assume that there is a free Hamiltonian and isometric action of $U (1) $ on $M$ generated by the  vector field $\xi $.  Let $\mu$ be the  associated Hamiltonian, $i _\xi  \omega  =- \mathrm{d} \mu $. Denote by $\Sigma _c $ the level set $\mu ^{-1} (c) $. 
It is well known that  if $c$ is a regular value of $\mu$ then the quotient $ \Sigma _c / U (1) $ equipped with the metric  and K\"ahler form induced by $g$ and $\omega$  is a smooth K\"ahler manifold of complex dimension $\operatorname{dim} _\mathbb{C}   M -1 $. 
The vector fields $  \xi  $, $J \xi  $ are respectively tangent and orthogonal to  $\Sigma _c $. In particular $ J \xi  $ spans the $g$-orthogonal complement of $T \Sigma _c $. Therefore the quotient metric $\overline{g }$ can be identified with the $g$-orthogonal complement of $(\xi  , J \xi  ) $ in $M$, 
\begin{equation}
\label{gdec} 
g = \bar g +  g ^\perp 
\end{equation} 
where $g ^\perp =g |_{ \operatorname{Span} (\xi , J\xi  )} $. Writing  $ \theta = g (\xi  , \cdot )$, $\theta _J =g ( J \xi  , \cdot )$ 
we have
\begin{equation}
g = \bar g    + \frac{\theta ^2 +  \theta _J ^2}{g (\xi  , \xi  )}.
\end{equation}

Take now  $ M =\mathbb{C}  ^2 $ and consider the $U (1) $ action is generated by the nowhere vanishing vector field $\xi =X _3 $,  which preserves both $\gc $ and $\oc $. 
The action is Hamiltonian with
\begin{equation}
\label{c2ham} 
\mu  =  \frac{1}{2} ( |z _1 |^2 +  |z _2 |^2 ).
\end{equation} 
Taking a basis $(X _1  , X _2 =J X _1  )$ for the orthogonal complement of $(X _3 , X _4 =J X _3 )$ as in (\ref{adaptedframe}),  gives back (\ref{gc2split}),
\begin{equation}
\begin{split} 
\gc&= \bar g +   g ^\perp  
= \frac{\theta  _1  ^2 + \theta  _2 ^2 }{|X _1 |^2 }  +  \frac{\theta _3 ^2  + \theta _4 ^2 }{ |X_3 |^2 } 
 =\frac{ |z _1 \mathrm{d} z _2-  z _2 \mathrm{d} z _1 |^2 }{|z _1 |^2 + |z _2 |^2 }  + \frac{  |\bar{z} _1  \mathrm{d} z _1 + \bar{z} _2 \mathrm{d} z _2 |^2 }{{|z _1 |^2 + |z _2 |^2 }} .
\end{split} 
\end{equation} 
%\begin{equation}
%\begin{split}
%\operatorname{Re} (\bar{z} _1  \mathrm{d} z _1 + \bar{z} _2 \mathrm{d} z _2) &
%= \operatorname{Re} (\bar\zeta \mathrm{d} \zeta ) = r \mathrm{d} r ,\\
%\operatorname{Im} (\bar{z} _1  \mathrm{d} z _1 + \bar{z} _2 \mathrm{d} z _2) &
%= \operatorname{Im} (\bar\zeta \mathrm{d} \zeta ) + \frac{|\zeta |^2 }{1 + |z| ^2 } \operatorname{Im} (\bar{z} \mathrm{d} z)
%= \frac{ r ^2 }{2}\eta _3,
%\end{split}
%\end{equation} 
For $R \neq 0 $ the $\mu$ level sets  $|z _1 |^2  + |z _2 |^2 = R ^2 $ are 3-spheres and the induced metric is
\begin{equation} 
\begin{split} 
g _{  S ^3 _R }&
= \frac{\theta  _1  ^2 + \theta  _2 ^2 }{|X _1 |^2 }  +  \frac{\theta _3 ^2  }{ |X_3 |^2 } 
 =\frac{ R^2}{4} \left[ \frac{ 4 |\mathrm{d} w|^2}{(1+|w|^2)^2} + 4\left(  \frac{\operatorname{Im} (\bar\zeta \mathrm{d} \zeta) }{R^2 } + a \right)  ^2 \right]
 = \frac{R ^2 }{4} (\eta _1 ^2 + \eta _2 ^2 + \eta _3 ^2 ),
 \end{split} 
\end{equation} 
having used (\ref{sphercoordsc4}).

\section{Calabi's construction}
\label{calabiconstr} 
For reference, we summarise here Calabi's construction \cite{calabi:1979} of a Ricci-flat K\"ahler manifold on the canonical bundle of a K\"ahler-Einstein manifold with scalar curvature $s \neq 0 $. We partially follow the presentation given in \cite{salamon:1990}.

\subsection{The Chern connection}
\label{chernconn} 
We recall, see e.g.~\cite{moroianu:2007}, that a holomorphic vector bundle  $E \rightarrow M$ equipped with a Hermitian product $\langle \cdot , \cdot \rangle $ carries a preferred connection $\nabla  $ known as the Chern connection.
It  is characterised by  compatibility with the Hermitian product: $\nabla _X  \langle \sigma _1 , \sigma _2 \rangle  = \langle X (\sigma _1 ), \sigma _2 \rangle + \langle \sigma _1 , X (\sigma _2 ) \rangle $ and compatibility with the complex structure on the total space.
% $E$: $\nabla _{ \bar Z } \sigma =0 $ for any anti-holomorphic vector field $\bar Z \in \Gamma (T ^{ 0,1 } M )$  and holomorphic section $\sigma\in \Gamma (E ^{ 1,0 } ) $.
%Let $r$ be the rank of $E$, $( \sigma _i )$ be a  local frame  and $(a ^i _j )$ the corresponding connection coefficients. Since the Chern connection is compatible with the complex structure it preserves the degree of objects it acts onto, so the $\mathfrak{gl } (r , \mathbb{C}  )$-valued local 1-forms are defined by
%\begin{equation}
%\nabla \sigma _j  =  a ^i _j \sigma _i, \quad \nabla \bar \sigma _{ \bar j}  =  \bar a ^{\bar i} _ {\bar j} \bar \sigma _{\bar i}
%\end{equation} 
%with no mixed coefficients.
%Written out in components, compatibility with the Hermitian product gives
%\begin{equation}
%\mathrm{d} ( h (\sigma _a , \sigma _b )  )=
%\partial h _{ a \bar b }+ \overline{\partial } h _{ a \bar b } 
%=h (\nabla \sigma _a , \sigma _b )   + h (\sigma _a , \nabla \sigma _b )  
%=  a ^c _a h _{ c \bar b }  +  \overline{a ^c _a h _{ c \bar b } },
%\end{equation} 
%where the second term  is the transpose conjugate of the first one. If $(\sigma _i ) $ is a holomorphic frame then $\nabla \sigma  _i = \nabla ^{ 1,0 } \sigma _i + \overline{\partial }\sigma _i = \nabla ^{ 1,0 } \sigma _i  $ so $a$ is of type $(1,0 )$ and matching the form degrees we get
%\begin{equation}
%a ^c _b =  h ^{ c \bar d} \partial h _{ b\bar d}.
%\end{equation} 
%If instead we are using a unitary frame $(e _i ) $ then $\mathrm{d} (h (\sigma _a , \sigma _b )  )=0 $ and so the connection coefficients are anti-hermitian.
If $E$ is a line bundle, $\sigma$ a section of $E$  then for some connection form $\a$
\begin{equation}
\nabla \sigma = \a\otimes  \sigma.
\end{equation} 
Compatibility with the Hermitian product implies
\begin{equation}
\label{chernline} 
 \mathrm{d}  \log \| \sigma \| ^2  = \a + \bar{\a }.
\end{equation} 
For a holomorphic section (\ref{chernline}) implies $\a =\partial \log \| \sigma \| ^2  $. If instead $\sigma =e $ is unitary, that is   $\|e \| ^2 =1 $, then  the connection is purely imaginary,  $\a =- \bar \a $. We recall that if $\a$ is the Chern connection on a line bundle $E$ then the $- \a $ is the Chern connection on $E ^\ast $.

Assume now that $M$ is  K\"ahler and  $E = T M \otimes \mathbb{C}$ is the complexified tangent bundle. Denote by $K =\Lambda  ^{ n,0 } (M) $ the canonical bundle of $M$, by $K ^\ast $ its dual and let $\a$ be the Chern connection on $K$, $\mathrm{d} \a $ the corresponding curvature.  Since $M$ is K\"ahler, the Chern connection on $E$ is just the Levi-Civita connection on $T M $ extended by complex linearity. The dual $K ^\ast $ of the canonical bundle $K =\Lambda  ^{ n,0 } (M) $  of $M$ is  the determinant bundle of $E $,  so  the  Chern connection  on $K ^\ast $ is  the trace of the Levi-Civita connection on $TM$, and the Chern connection $\a$ on $K$ its opposite. Therefore
\begin{equation}
\mathrm{d} \a =i \rho _M 
\end{equation} 
where $\rho _M $ is the Ricci form on $M$.
  If  $M$ is K\"ahler-Einstein of complex dimension $n$ and scalar curvature $s$ then $\rho _M = \tfrac{s}{2n} \omega _M  $, where $\omega _M $ is the K\"ahler  form  on $M$, hence
\begin{equation} 
- i \mathrm{d} \a =\rho _M = \frac{s}{2n} \omega _M .
\end{equation} 

\subsection{Calabi's construction}
\label{calabiconstr} 
In this section $M$ is a K\"ahler-Einstein  manifold of complex dimension $n$ with Riemannian metric $g _M $, K\"ahler form $\omega _M $, Riemannian volume element 
\begin{equation} \mathrm{vol} _M  =\frac{\omega ^n }{n!}
\end{equation}
and scalar curvature $s$. The bundle $K$ is the canonical bundle of $M$, $\a$ the connection form of the Chern connection $\nabla $ on $K$. Compatibility of $\nabla $ with the complex structure is equivalent to  $\nabla ^{ 0,1 } =\overline{\partial }$.

Let $\sigma$ be a unitary section of $K$. Since $\sigma$ is a form of degree $(n,0 )$, $\partial \sigma =0  = \nabla  ^{ 1,0 } \sigma$ and so
\begin{equation}
\label{de} 
\mathrm{d} \sigma = \overline{\partial } \sigma = \nabla \sigma = \a \wedge \sigma.
\end{equation} 
The Riemannian metric on $M$ induces a Hermitian product $\langle \cdot , \cdot \rangle $ on $K$ given by
\begin{equation}
\label{fibreproduct} 
\sigma  \wedge * \overline{ \sigma  } = \| \sigma \| ^2  \mathrm{vol} _M .
\end{equation} 
Since  $* : \Lambda ^{ p,q }(M) \rightarrow \Lambda ^{ n-q,n-p }(M) $ and $* ^2 =(-1 )^{  (p + q) ^2 } $,  $* $ acts on  $\Lambda ^{ 0,n }(M)   $ by multiplication by $ i  ^{n ^2 } $. Hence (\ref{fibreproduct})  is equivalent to
\begin{equation}
\label{fibreproduct2} 
\sigma  \wedge \overline{\sigma  } = (-i )^{ n ^2 }  \|\sigma \| ^2  \frac{  \omega _M  ^{n} }{n!}.
\end{equation} 

%Since $M$ is K\"ahler-Einstein, its Ricci form satisfies $ \rho = \tfrac{s}{2n}  \omega _M  $. The Chern connection is the trace of the Levi-Civita one, so the Chern curvature $\mathrm{d} a $ satisfies
%\begin{equation}
%i \mathrm{d} a = \frac{s}{2n} \omega  _M .
%\end{equation} 

Let $\zeta$ be a coordinate on the fibres of $K$, set
\begin{equation}
\theta =\mathrm{d} \zeta +  \zeta \a.
\end{equation} 
and define 
\begin{align}
\label{tildeomega} 
 \omega &= \lambda \left[ u\,  \omega _M  +   \frac{2in }{s} u ^\prime \theta \wedge \bar \theta \right] ,\\
 \label{Omegadef} 
\Omega  &= \theta \wedge e ,
\end{align} 
where $u $ is a function of $|\zeta |^2 $ and $\lambda$ is some non-zero constant. 
%\begin{equation} 
%D =\left( \frac{s (1 + n) }{2n} \right) ^{ 1/(n + 1 ) }
%\end{equation} 
We are now going to show that $(\omega , \Lambda )$ satisfy
\begin{align}
\label{closomega} 
\mathrm{d} \omega &=0, \\
\label{closOmega} 
 \mathrm{d}  \Omega  &=0, \\
 \label{wedge0} 
 \Omega  \wedge \omega &=0,\\
 \label{omeganorm} 
\Omega  \wedge \bar \Omega  &=C \omega ^{n + 1 }.
\end{align} 
Closure of $\Omega$  means that the complex structure defined by declaring $\theta$ to be a $(1,0)$ form on $K$ is integrable, while (\ref{wedge0})  implies that ${\omega}$ is a $(1,1)$ form in the complex structure defined by $\Omega$. Equation (\ref{closomega}) and non-degeneracy of $\omega$ then imply that $\omega$ is a K\"ahler form on $K$. Finally the normalisation condition  (\ref{omeganorm}), where $C$ is some constant, implies that the K\"ahler metric on $K$ is Ricci-flat.

Since (\ref{omeganorm}) is equivalent to
\begin{equation}
\| \Omega  \| ^2 =i ^{ (n  + 1 )^2} C  ,
\end{equation} 
the value of the  constant $C $ depends on the value $\| \Omega  \| $. Here we follow \cite{hitchin:1997} by taking
\begin{equation} 
\label{convC} 
C =(-i )^{ (n + 1 )^2 } 
= \begin{cases}
 1 \quad  & \text{if $n$ is odd},\\
  -i \quad  & \text{if $n$ is even},
  \end{cases} 
  \end{equation} 
corresponding to  $\| \Omega  \| ^2 =1  $.

We now show the properties (\ref{closomega})--(\ref{omeganorm}) hold for $\omega$, $\Omega$ defined by (\ref{tildeomega}), (\ref{Omegadef}) and $u$ given below by (\ref{usol}). One calculates
\begin{equation}
\begin{split} 
\label{dzeta2} 
\mathrm{d} |\zeta |^2 &= \zeta \bar \theta + \bar \zeta \theta, \\
\mathrm{d} \theta &= \zeta \mathrm{d} \a +  \mathrm{d} \zeta \wedge \a  
= i \frac{s}{2n} \zeta \omega_M +  \mathrm{d} \zeta \wedge \a 
= i \frac{s}{2n} \zeta \omega_M  +  \theta  \wedge \a ,
\end{split} 
\end{equation} 
from which  $ \mathrm{d}  \omega =0 $ easily follows.
Since, using (\ref{de}),
\begin{equation}
  \Omega  = (   \mathrm{d} \zeta +  \zeta \a )\wedge e =\mathrm{d} (\zeta e )
\end{equation} 
$\Omega  $ is closed. Since $\Omega  $   is  of type $(n + 1,0 )$ one has $\omega _M \wedge \Omega  =0 $ and $ \omega \wedge \Omega  =0 $ follows.  
Using (\ref{fibreproduct2}) applied to the K\"ahler metric determined by $ \omega $ one has
\begin{equation}
\Omega  \wedge \bar \Omega   = (-i )^{( n + 1 )^2}  \| \Omega  \| ^2 \frac{  \omega ^{ n + 1 }}{ (n+1)!},
\end{equation} 
and from (\ref{tildeomega}) 
\begin{equation}
  \omega ^{ n + 1 }  
= \frac{2i n }{s}\lambda ^{ n + 1 } u ^n u ^\prime \omega _M ^n \wedge \theta \wedge \bar \theta 
%  =  i (n + 1 ) ^2 u ^n u ^\prime   \omega _M^n \wedge \theta \wedge \bar \theta ,
\end{equation} 
giving
\begin{equation}
\Omega  \wedge \bar \Omega  
= (-1 )^n (-i )^{ n ^2 }\frac{2n}{s (n + 1 )} \lambda ^{ n + 1 }  \| \Omega  \| ^2u ^n u ^\prime \mathrm{vol} _M \wedge \theta \wedge \bar \theta 
%= (-1 )^{n}  (-i )^{ n ^2 }   \| \Lambda \| ^2 (n + 1 ) u ^n u ^\prime \mathrm{vol} _M \wedge \theta \wedge \bar \theta .
\end{equation} 
We can also compute directly
\begin{equation}
\begin{split} 
\Omega  \wedge \bar \Omega  &
= (-1 )^n \theta \wedge \bar \theta \wedge e \wedge \bar e
=(-1 )^n (-i )^{ n ^2 }\theta \wedge \bar \theta \wedge e \wedge *  \bar e \\ &
=(-1 )^n (-i )^{ n ^2 } \mathrm{vol} _M \wedge \theta \wedge \bar \theta .
\end{split} 
\end{equation} 
Thus
\begin{equation}
\frac{2n}{s (n + 1 )} \lambda ^{ n + 1 }  \| \Omega  \| ^2  u ^n u ^\prime =1
%(n + 1 )\| \Lambda \| ^2 u ^n u ^\prime = 1
\end{equation} 
and in order to satisfy the normalisation condition $\| \Omega  \| ^2 =1 $  we need to impose
\begin{equation}
\label{reluprime} 
\frac{2n}{s (n + 1 )} \lambda ^{ n + 1 }  u ^n u ^\prime  =  1, 
\end{equation} 
which has solution
\begin{equation}
\label{usol} 
u =  \left(  c | \zeta |^2  +  \kappa  \right)^{ \frac{1}{n + 1} },
\end{equation} 
where $\kappa $ is an integration constant and
\begin{equation}
c =\frac{s (n + 1 )^2 }{2n \lambda ^{ n + 1 }}.
\end{equation} 

In conclusion   the K\"ahler form $\omega$ and metric $g$  on $K$ are
\begin{align}
 \omega &=\lambda u\,  \omega_M +i (n+1) (\lambda u )^{ -n } \theta \wedge \bar \theta ,\\
g &
= \lambda u \, g _M  + 2(n+1)(\lambda u )^{ -n } |\theta  |^2 .
\end{align} 
If $s>0 $, $\kappa >0 $ and $g _M $ is complete then $g$ is globally defined on the total space of $K$ and complete.

\newpage 
\bibliographystyle{amsplain}
%\bibliography{/Users/snoopy/Dropbox/Academia/Papers/biblio.bib}

\providecommand{\bysame}{\leavevmode\hbox to3em{\hrulefill}\thinspace}
\providecommand{\MR}{\relax\ifhmode\unskip\space\fi MR }
\providecommand{\MRhref}[2]{%
  \href{http://www.ams.org/mathscinet-getitem?mr=#1}{#2}
}
\providecommand{\href}[2]{#2}

\end{document}